\documentclass{article}

\usepackage{amssymb,amsmath,amscd,amsthm,epsfig}

\newcommand{\bbq}{{\mathbf Q}}

\newcommand{\inv}{\mathrm{inv}}
\newcommand{\Des}{\mathrm{Des}}

\newcommand{\maj}{\mathrm{maj}}
\newcommand{\rmaj}{\mathrm{comaj}}

\newcommand{\dm}{\mathrm{dim\ }}
\newcommand{\kr}{\mathrm{ker\ }}
\newcommand{\supp}{\mathrm{Supp}}

\newcommand{\DP}{{\cal P}}

\newcommand{\sn}{\mathrm{span}}
\newcommand{\la}{{\lambda}}
\newcommand{\then}{{\Longrightarrow}}

\def \ox {{\overline x}}
\def \oy {{\overline y}}
\def \oz {{\overline z}}

\newtheorem{thm}{Theorem}[section]
\newtheorem{pro}[thm]{Proposition}
\newtheorem{lem}[thm]{Lemma}
\newtheorem{cla}[thm]{Claim}

\newtheorem{cor}[thm]{Corollary}

\newtheorem{obs}[thm]{Observation}

\newtheorem{exa}[thm]{Example}
\newtheorem{df}[thm]{Definition}
\newtheorem{rem}[thm]{Remark}

\newenvironment{note}[1]{\par\addvspace{\medskipamount}\noindent
                         {\bf {#1}}\sl
                       }{\par\addvspace{\medskipamount}\rm}

\newcommand{\beq}{\begin{equation}}
\newcommand{\eeq}{\end{equation}}
\newcommand{\beqn}{\begin{eqnarray}}
\newcommand{\eeqn}{\end{eqnarray}}

\newcommand{\htop}{\hbox{\rm top}}

\def \ox {{\overline x}}
\def \la {\lambda}

\title{The Combinatorics of the Garsia-Haiman Modules\\ for Hook
Shapes}
\author{%
Ron M.\ Adin\thanks{%
        Supported in part by the Israel Science Foundation,
        grant no.\ 947/04.}\\
        Department of Mathematics\\
        Bar-Ilan University\\ Ramat-Gan 52900, Israel\\
        radin@math.biu.ac.il\and
Jeffrey B.\ Remmel\thanks{%
        Supported in part by NSF grant DMS 0400507.}\\
        Department of Mathematics\\
        University of California, San Diego\\
        La Jolla, CA 92093\\
        remmel@math.ucsd.edu\and
Yuval Roichman\thanks{%
        Supported in part by the Israel Science Foundation,
        grant no.\ 947/04,
        and by the University of California, San Diego.}\\
        Department of Mathematics\\
        Bar-Ilan University\\ Ramat-Gan 52900, Israel\\
        yuvalr@math.biu.ac.il}
\date{February 16, 2007}

\begin{document}
\maketitle

\begin{abstract}
Several bases of the Garsia-Haiman modules for hook shapes are
given, as well as combinatorial decomposition rules for these
modules. These bases and rules extend the classical ones for the
coinvariant algebra of type $A$. We also give a
decomposition of the Garsia-Haiman modules into descent representations.
\end{abstract}

\section{Introduction}

\subsection{Outline}

In~\cite{GH93}, Garsia and Haiman introduced a module ${\bf H}_\mu$ for each
partition $\mu$, which
we shall call the Garsia-Haiman module for $\mu$. Garsia and Haiman
introduced the modules ${\bf H}_\mu$ in
attempt to prove Macdonald's $q,t$-Kostka polynomial conjecture and,
in fact, the modules ${\bf H}_\mu$ played a major role in the
resolution of Macdonald's conjecture \cite{Hai01}.
When the shape $\mu$ has a single row, this module is isomorphic
to the coinvariant algebra of type $A$. Our goal here is to
understand the structure of this module when $\mu$ is a hook shape
$(1^{k-1},n-k+1)$.

A family of bases for the  Garsia-Haiman module of hook shape
$(1^{k-1},n-k+1)$
is presented. This family includes the {\em $k$-th Artin basis}, the
{\em $k$-th descent basis}, the {\em $k$-th Haglund basis} and the
{\em $k$-th Schubert basis} as well as other bases. While the
first basis appears in~\cite{GH1}, the others are new and have
interesting applications.

The $k$-th Haglund basis realizes Haglund's statistics for the modified
Macdonald polynomials in the hook
case. The $k$-th descent basis extends the well known Garsia-Stanton
descent basis for the coinvariant algebra. The advantage of the
$k$-th descent basis is that the $S_n$-action on it may be
described explicitly. This description implies combinatorial rules
for decomposing the bi-graded components of the module into
Solomon descent representations and into irreducibles. In
particular, a constructive proof of a formula due to Stembridge is
deduced.

\medskip

\tableofcontents

%

\subsection{General Background}\label{background}

In 1988, I.\ G.\ Macdonald~\cite{Mac88} introduced a remarkable
new basis for the space of  symmetric functions. The elements of
this basis are denoted $P_{\lambda}(\ox;q,t)$, where $\la$ is a
partition, $\ox$ is a vector of indeterminates, and $q,t$ are
parameters. The $P_{\lambda}(\ox;q,t)$'s,  which are now called
``{\em Macdonald polynomials}'', specialize to many of the
well-known bases for the symmetric functions, by suitable substitutions
for the parameters $q$ and $t$. In fact, we can obtain in this
manner the Schur functions, the Hall-Littlewood symmetric
functions, the Jack symmetric functions, the zonal symmetric
functions, the zonal spherical functions, and the elementary and
monomial symmetric functions.

\begin{figure}[ht]
\begin{center}
\epsfig{file=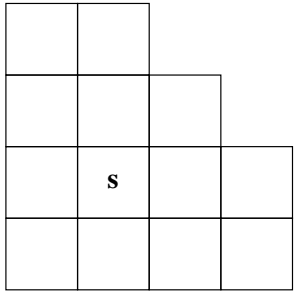}
\caption{Diagram of a partition.}
\label{fig:Mac1}
\end{center}
\end{figure}

\medskip

Given a cell $s$ in the Young  diagram (drawn according to the
French convention) of a partition $\la$, let $leg_\la(s)$,
$leg'_\la(s)$, $arm_\la(s)$, and $arm'_\la(s)$ denote the number
of squares that lie above, below, to the right, and to left of $s$
in $\la$, respectively. For example, when $\la =(2,3,3,4)$ and $s$
is the cell pictured in Figure \ref{fig:Mac1}, $leg_\la(s)= 2$,
$leg'_\la(s)= 1$, $arm_\la(s)= 2$ and $arm'_\la(s)= 1$. For each
partition $\la$, define
\begin{equation}
h_\la(q,t) := \prod_{s \in \la} (1 - q^{arm_\la(s)}
t^{leg_\la(s)+1})
\end{equation}
For a partition $\la = (\la_1,\ldots,\la_k)$ where $0<\la_1 \leq
\ldots \leq \la_k$, let $n(\la) := \sum_{i=1}^k (k-i)\la_i$.
Macdonald introduced the $(q,t)$-Kostka polynomials
$K_{\la,\mu}(q,t)$ via the equation
\begin{equation}\label{qtKost1}
J_\mu(\ox;q,t) = h_\mu(q,t)P_\mu(\ox;q,t) = \sum_\la
K_{\la,\mu}(q,t) s_\la[X(1-t)],
\end{equation}
and conjectured that they are polynomials in $q$ and $t$ with
non-negative integer coefficients.

\begin{figure}[ht]
\begin{center}
\epsfig{file=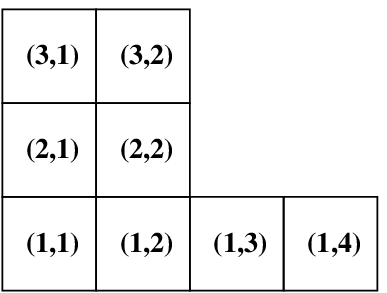}
\caption{Labelling of the cells of a partition.}
\label{fig:Mac2}
\end{center}
\end{figure}

In an attempt to prove Macdonald's conjecture, Garsia and Haiman
\cite{GH93} introduced the so-called {\em modified Macdonanld
polynomials} $\tilde{H}_\mu(\ox;q,t)$ as
\begin{equation}\label{modMac}
\tilde{H}_\mu(\ox;q,t) = \sum_\la \tilde{K}_{\la,\mu}(q,t)
s_\la(\ox),
\end{equation}
where $\tilde{K}_{\la,\mu}(q,t) :=t^{n(\mu)}K_{\la,\mu}(q,1/t)$.
Their idea was that $\tilde{H}_\mu(\ox;q,t)$ is the Frobenius
image of the character generating function of a certain bi-graded
module ${\bf H}_\mu$ under the diagonal action of the symmetric
group $S_n$. To define ${\bf H}_\mu$, assign
$(row,column)$-coordinates to squares in the first quadrant,
obtained by permuting the $(x,y)$ coordinates of the upper
right-hand corner of the square so that the lower left-hand square
has coordinates (1,1), the square above it has coordinates (2,1),
the square to its right has coordinates (1,2), etc. The first
(row) coordinate of a square $w$ is denoted $row(w)$, and the
second (column) coordinate of $w$ is the denoted $col(w)$. Given a
partition $\mu \vdash n$, let $\mu$ also denote the corresponding
Young diagram, drawn according to the French convention, which
consists of all the squares with coordinates $(i,j)$ such that $1
\leq i \leq \ell(\mu)$ and $1 \leq j \leq \mu_i$. For example, for
$\mu = (2,2,4)$, the labelling of squares is depicted in Figure
\ref{fig:Mac2}.

Fix an ordering $w_1, \ldots, w_n$ of the squares of $\mu$, and
let
\begin{equation}\label{Dmu}
\Delta_\mu(x_1, \ldots, x_n; y_1, \ldots, y_n): = \det
\left(x_i^{row(w_j)-1}y_i^{col(w_j)-1}\right)_{i,j}.
\end{equation}
For example,
$$
\Delta_{(2,2,4)}(x_1, \ldots, x_8; y_1, \ldots, y_8) = \det \left(
\begin{matrix} 1 & y_1 & y_1^2 & y_1^3 &x_1 & x_1 y_1 & x_1^{2} & x_1^{2}y_1 \\
1 & y_2 & y_2^2 & y_2^3 &x_2 & x_2 y_2 & x_2^{2} & x_2^{2}y_2 \\
\vdots & & & & & & &\vdots\\
1 & y_8 & y_8^2 & y_8^3  & x_8 & x_8 y_8 & x_8^{2} & x_8^{2} y_8
\end{matrix} \right)
$$

Now let ${\bf H}_\mu$ be the vector space of polynomials spanned
by all the partial derivatives of \\
$\Delta_\mu(x_1, \ldots,x_n;y_1,
\ldots, y_n)$. The symmetric group $S_n$ acts on ${\bf H}_\mu$
diagonally, where for any polynomial $P(x_1, \ldots, x_n; y_1,
\ldots, y_n)$ and any permutation $\sigma \in S_n$,
$$
P(x_1, \ldots, x_n; y_1, \ldots, y_n)^{\sigma} := P(x_{\sigma_1},
\ldots, x_{\sigma_n}; y_{\sigma_1}, \ldots, y_{\sigma_n}).
$$
The bi-degree $(h,k)$ of a monomial $x_1^{p_1} \cdots x_n^{p_n}
y_1^{q_1} \cdots y_n^{q_n}$ is defined by $h := \sum_{i=1}^n p_i$
and $k := \sum_{i=1}^n q_i$. Let ${\bf H}_\mu^{(h,k)}$ denote
space of homogeneous polynomials of degree $(h,k)$ in ${\bf
H}_\mu$. Then
$$
{\bf H}_\mu = \bigoplus_{(h,k)} {\bf H}_\mu^{(h,k)}.
$$
The $S_n$-action clearly preserves the bi-degree so that $S_n$
acts on each homogeneous component ${\bf H}_\mu^{(h,k)}$. The
character of the $S_n$-action on ${\bf H}_\mu^{(h,k)}$ can be
decomposed as
\begin{equation}\label{char1}
\chi^{{\bf H}_\mu^{(h,k)}} = \sum_{\la \vdash n}
\chi_{\la,\mu}^{(h,k)} \chi^{\la},
\end{equation}
where $\chi^\la$ is the irreducible character of $S_n$ indexed by
the partition $\la$ and the  $\chi_{\la,\mu}^{(h,k)}$'s are
non-negative integers. We then define the character generating
function of ${\bf H}_\mu$ to be
\begin{eqnarray}\label{char2}
\chi^{{\bf H}_\mu}(q,t)
&=& \sum_{h,k \geq 0} q^h t^k
    \sum_{\la \vdash |\mu|} \chi_{\la,\mu}^{(h,k)} \chi^\la \\
&=& \sum_{\la \vdash |\mu|} \chi^\la
    \sum_{h,k \geq 0} \chi_{\la,\mu}^{(h,k)} q^h t^k \nonumber
\end{eqnarray}
The Frobenius map $F$, which maps the center of the group algebra
of $S_n$ to $\Lambda_n(\ox)$, is defined by sending the character
$\chi^\la$ to the Schur function $s_\la(\ox)$. Garsia and Haiman
conjectured that the Frobenius image of $\chi^{{\bf H}_\mu}(q,t)$
which they denoted by $F_{\mu}(q,t)$ is the modified Macdonald
polynomial $\tilde{H}_{\mu}(\ox;q,t)$. That is, they conjectured
that
\begin{eqnarray}\label{char3}
F(\chi^{{\bf H}_{\mu}}(q,t))
&=& \sum_{h,k \geq 0} q^h t^k
    \sum_{\la \vdash |\mu|} \chi_{\la,\mu}^{(h,k)} s_\la(\ox) \\
&=& \sum_{\la \vdash |\mu|} s_\la(\ox) \sum_{h,k \geq 0}
    \chi_{\la,\mu}^{(h,k)} q^h t^k \nonumber \\
&=& \tilde{H}_{\mu}(\ox;q,t) \nonumber
\end{eqnarray}

so that
\begin{equation}\label{char4}
\tilde{K}_{\la,\mu}(q,t) = \sum_{h,k \geq 0} \chi_{\la,\mu}^{(h,k)} q^h t^k
\end{equation}
Since Macdonald proved that $K_{\la,\mu}(1,1)= f_\la$,
the number of standard tableau of shape $\la$,
equations (\ref{char3}) and (\ref{char4}) led Garsia and Haiman
\cite{GH93} to conjecture that  as an $S_n$-module, ${\bf H}_\mu$
carries the regular representation.
This conjecture was eventually proved by Haiman \cite{Hai01}
using the algebraic geometry of the Hilbert Scheme.

The goal of this paper is to understand the structure of the
modules ${\bf H}_\mu$ when $\mu$ is a hook shape
$(1^{k-1},n-k+1)$. These modules were studied before by
Stembridge~\cite{Ste94}, Garsia and Haiman~\cite{GH1},
Allen~\cite{Allen2} and Aval~\cite{Aval}. This paper suggests a
detailed combinatorial analysis of the modules.

\subsection{Main Results - Bases}\label{s.main}


Consider the inner product $\langle \ ,\ \rangle$ on the
polynomial ring $\bbq[\bar x,\bar
y]=\bbq[x_1,\dots,x_n,y_1,\dots,y_n]$ defined as follows:
\begin{equation}
\langle f,g \rangle := \hbox{\rm constant term of }
f(\partial_{x_1},\dots,\partial_{x_n};\partial_{y_1},\dots,\partial_{y_n})
g(x_1,\dots,x_n;y_1,\dots,y_n)\qquad(\forall f,g\in Q_n),
\end{equation}
where $f(\partial_{x_1},\dots,\partial{y_n})$ is the differential operator
obtained by replacing each variable $x_i$ ($y_i$) in $f$ by the corresponding
partial derivative $\frac{\partial}{\partial {x_i}}$
($\frac{\partial}{\partial{y_i}}$).
Let  $J_\mu$ be the $S_n$-module dual to ${\bf H}_\mu$ with
respect to $\langle\ ,\ \rangle$, and let ${\bf H}_\mu':=\bbq[\bar
x,\bar y]/J_\mu$. It is not difficult to see that ${\bf H}_\mu$
and ${\bf H}_\mu'$ are isomorphic as $S_n$-modules.

\subsubsection{The $k$-th Descent Basis}\label{k-des}

The {\em descent set} of a permutation $\pi\in S_n$ is
$$
\Des(\pi):=\{i\,:\,\pi(i)>\pi(i+1)\}.
$$
Garsia and Stanton~\cite{GS} associated with each $\pi\in S_n$, the
{\em descent monomial}
$$
a_\pi := \prod_{i\in\Des(\pi)} (x_{\pi(1)} \cdots x_{\pi(i)}) =
\prod\limits_{j=1}^{n-1} x_{\pi(j)}^{|\Des(\pi)\cap
\{j,\ldots,n-1\}|}.
$$
Using Stanley-Reisner rings, Garsia and Stanton~\cite{GS} showed
that the set $\{a_\pi\,:\,\pi\in S_n\}$ forms a basis for the
coinvariant algebra of type $A$. See also~\cite{Stg}
and~\cite{Allen1}.

\medskip

\begin{df}
For every integer $1\leq k\leq n$ and permutation $\pi\in S_n$
define
$$
d_i^{(k)}(\pi) := \begin{cases} |\Des(\pi)\cap \{i,\ldots,k-1\}|,
&\hbox{if $1\leq i<k$};\cr 0, &\hbox{if $i=k$};\cr |\Des(\pi)\cap
\{k,\ldots,i-1\}|, &\hbox{if $k< i\leq n$}.
\end{cases}
$$
\end{df}

\begin{df}\label{d.basis}
For every integer $1\leq k\leq n$ and permutation $\pi\in S_n$
define the {\em $k$-th descent monomial}
\begin{eqnarray*}
a_\pi^{(k)} &:=& \prod_{i\in\Des(\pi) \atop i\leq k-1} (x_{\pi(1)}
\cdots x_{\pi(i)}) \cdot
\prod_{i\in\Des(\pi) \atop i\ge k} (y_{\pi(i+1)} \cdots y_{\pi(n)})\\
&=& \prod\limits_{i=1}^{k-1} x_{\pi(i)}^{d_i^{(k)}(\pi)} \cdot
\prod\limits_{i=k+1}^{n} y_{\pi(i)}^{d_i^{(k)}(\pi)}.
\end{eqnarray*}
\end{df}

For example, if $n=8$, $k=4$, and $\pi = 8~6~1~4~7~3~5~2$, then
$Des(\pi) = \{1,2,5,7\}$,
$(d_1^{(4)}(\pi), \ldots, d_8^{(4)}(\pi)) = (2,1,0,0,0,1,1,2)$,
and $a^{(4)}_{\pi} = x_1^2 x_2 y_6 y_7 y_8^2$.

Note that $a_\pi^{(n)}=a_\pi$, the Garsia-Stanton descent
monomial.

Consider the partition $\mu = (1^{k-1},n-k+1)$.

\begin{thm}\label{m1}
For every $1\leq k\leq n$, the set of $k$-th descent monomials
$\{a_\pi^{(k)}\,:\, \pi\in S_n\}$ forms a basis for the
Garsia-Haiman module ${\bf H}_{(1^{k-1},n-k+1)}$.
\end{thm}

Two proofs of Theorem~\ref{m1} are given in this paper. In
Section~\ref{s.straight} it is proved
via a straightening algorithm. This proof implies an explicit
description of the Garsia-Haiman hook module ${\bf H}_{(1^{k-1},n-k+1)}'$.



\begin{thm}\label{t.m2}
For $\mu=(1^{k-1},n-k+1)$ the ideal $J_\mu={\bf H}_\mu^\perp$ defined
above is the ideal of $\bbq[\ox,\oy]$ generated by
\begin{itemize}
\item[(i)]
$\Lambda [\bar x]^+$ and $\Lambda [\bar y]^+$ (the symmetric functions
in $\bar x$ and $\bar y$ without a constant term),
\item[(ii)]
the monomials $x_{i_1}\cdots x_{i_k}$ ($i_1<\cdots<i_k$) and
$y_{i_1}\cdots y_{i_{n-k+1}}$ ($i_1<\cdots<i_{n-k+1}$), and
\item[(iii)]
the monomials $x_iy_i$ ($1\le i\le n$).
\end{itemize}
\end{thm}

This result has been obtained in a different form by J.-C.\
Aval~\cite[Theorem 2]{Aval}.

\subsubsection{The $k$-th Artin and Haglund Bases}\label{s.k-AH}

A second proof of Theorem~\ref{m1} is given in
Section~\ref{s.kick}. This proof applies a generalized version of
the Garsia-Haiman kicking process. This construction is extended
to a rich family of bases.

For every positive integer $n$, denote $[n]:=\{1,\dots,n\}$.
For every subset $A=\{i_1,\dots,i_k\}\subseteq [n]$ denote
$\bar x_A:=x_{i_1},\dots,x_{i_k}$ and $\bar y_A:=y_{i_1},\dots,y_{i_k}$.
Denote $\bar x:= \bar x_{[n]}=x_1,\dots,x_n$ and
$\bar y:= \bar y_{[n]}= y_1,\dots,y_n$.

Let $k,c\in [n]$, let $A = \{a_1,\dots,a_{k-1}\}$ be a subset
of size $k-1$ of $[n]\setminus c$, and let $\bar A:=[n]\setminus (A\cup\{c\})$.
Let $B_A$ be an arbitrary basis of the coinvariant algebra of $S_{k-1}$
acting on $\bbq[\bar x_A]$,
and let $C_{\bar A}$ be a basis of the coinvariant algebra of $S_{n-k}$
acting on $\bbq[\bar y_{\bar A}]$.
Finally define
$$
m_{(A,c, \bar A)}:=\prod\limits_{\{i\in A\, \, :\, i>c\}} x_i
\prod\limits_{\{j\in \bar A\, :\, j<c\}} y_j \in \bbq[\bar x,\bar
y].
$$
Then

\begin{thm}\label{mk}
The set
$$
\bigcup_{A,c}\ m_{(A,c, \bar A)} \ B_A\ C_{\bar
A}:=\bigcup_{A,c}\{ m_{(A,c, \bar A)}bc\ :\ b\in B_A, c\in C_{\bar
A}\}
$$
forms a basis for the Garsia-Haiman module ${\bf
H}_{(1^{k-1},n-k+1)}'$.
\end{thm}


\begin{df}
For every integer $1\leq k\leq n$ and permutation $\pi\in S_n$
define
$$
\inv_i^{(k)}(\pi) := \begin{cases} |\{j\,:\,i<j\le k \mbox{ and }
\pi(i)>\pi(j)\}|, &\hbox{if $1\le i<k$};\cr 0, &\hbox{if
$i=k$};\cr |\{j\,:\,k\leq j<i \mbox{ and } \pi(j)>\pi(i)\}|,
&\hbox{if $k< i\leq n$}.
\end{cases}
$$
%
For every integer $1\leq k\leq n$ and permutation $\pi\in S_n$
define the {\em $k$-th Artin monomial}
\begin{eqnarray*}
b_\pi^{(k)} &:=& \prod\limits_{i=1}^{k-1}
x_{\pi(i)}^{\inv_i^{(k)}(\pi)} \cdot \prod\limits_{i=k+1}^{n}
y_{\pi(i)}^{\inv_i^{(k)}(\pi)}.
\end{eqnarray*}
and the  {\em $k$-th Haglund monomial}
\begin{eqnarray*}
c_\pi^{(k)} &:=& \prod\limits_{i=1}^{k-1}
x_{\pi(i)}^{d_i^{(k)}(\pi)} \cdot \prod\limits_{i=k+1}^{n}
y_{\pi(i)}^{\inv_i^{(k)}(\pi)}.
\end{eqnarray*}
\end{df}

For example, if $n=8$, $k=4$, and $\pi = 8~6~1~4~7~3~5~2$, then
$Des(\pi) = \{1,2,5,7\}$,
$(inv_1^{(4)}(\pi), \ldots, inv_8^{(4)}(\pi)) = (3,2,0,0,1,0,2,0)$,
$b^{(4)}_{\pi} = x_1^3 x_2^2 y_5 y_7^2$, and
$c^{(4)}_{\pi} = x_1^2 x_2 y_5 y_7^2$.

Interesting special cases of Theorem~\ref{mk} are the following.

\begin{cor}\label{m11}
Each of the following sets :
$$
\{a_\pi^{(k)}\,:\,\pi\in S_n\},\ \ \ \
\{b_\pi^{(k)}\,:\,\pi\in S_n\},\ \ \ \
\{c_\pi^{(k)}\,:\,\pi\in S_n\}
$$
forms a basis for the Garsia-Haiman module ${\bf H}_{(1^{k-1},n-k+1)}'$.
\end{cor}

\begin{rem}\rm \
\begin{itemize}
\item[1.] Garsia and Haiman [12] showed that
$\{b_\pi^{(k)}\,:\,\pi\in S_n\}$ is a basis for ${\bf H}_{(1^{k-1},n-k+1)}'$.
Other bases of ${\bf H}_{(1^{k-1},n-k+1)}$ were also constructed
by J-C Aval \cite{Aval} and E. Allen \cite{Allen1,Allen2}. They used
completely different methods. Aval constructed a basis of the
form of an explicitly described set of partial differential operators
applied to $\Delta_{(1^{k-1},n-k+1)}$ and Allen constructed a basis for
${\bf H}_{(1^{k-1},n-k+1)}$ out his theory of bitableaux.

\item[2.] It should be noted that the last basis corresponds to
Haglund's statistics for the Hilbert series of
${\bf H}_{(1^{k-1},n-k+1)}$ that is implied by his combinatorial
interpretation for the modified Macdonald polynomial
$\tilde{H}_{(1^{k-1},n-k+1)}(\bar{x};q,t)$; see Section~\ref{QS}
below.
\item[3.] Choosing $B_A$ and $C_{\bar A}$ in
Theorem~\ref{mk} to be the Schubert bases of the coinvariant
algebras of $S_{k-1}$ (acting on $\bbq[\bar x_A]$) and of
$S_{n-k}$ (acting on $\bbq[\bar y_{\bar A}]$), respectively, gives
the {\em $k$-th Schubert basis}.
One may study the Hecke algebra actions on this basis along the
lines drawn in~\cite{APR1}.
\end{itemize}
\end{rem}

\subsection{Main Results - Representations}\label{rep}

\subsubsection{Decomposition into Descent
Representations}\label{m.des-rep}



The set of elements in a Coxeter group having a fixed descent set
carries a natural representation of the group, called a descent
representation. Descent representations of Weyl groups were first
introduced by Solomon \cite{So} as alternating sums of permutation
representations. This concept was extended to arbitrary Coxeter
groups, using a different construction, by Kazhdan and Lusztig
\cite{KL} \cite[\S 7.15]{Hum}. For Weyl groups of type $A$, these
representations also appear in the top homology of certain
(Cohen-Macaulay) rank-selected posets \cite{St82}. Another
description (for type $A$) is by means of zig-zag diagrams
\cite{Ge2, GR}. A new construction of descent representations for
Weyl groups of type $A$, using the coinvariant algebra as a
representation space, is given in~\cite{TAMS}.

\medskip

\begin{figure}[tbp]
\centering
$$
\epsfysize=.7in \epsffile{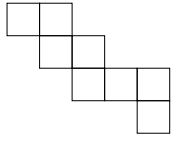}
$$
\caption{The zigzag shape corresponding to $n = 8$ and
$A = \{2,4,7\}$.}\label{fig:DA}\end{figure}

For every subset $A\subseteq \{1,\dots,n-1\}$, let
$$
S_n^A:=\{\pi\in S_n\, :\, \Des(\pi)=A\}
$$
be the corresponding {\em descent class} and let $\rho^A$ denote the
corresponding {\em descent representation} of $S_n$. Given
$n$ and subset $A = \{a_1 < \cdots < a_k\} \subseteq \{1,\dots,n-1\}$,
we can associate a composition of $n$, $comp(A) = (c_1,
\ldots, c_{k+1}) = (a_1, a_2 -a_1,
\ldots, a_k -a_{k-1},n-a_k)$ and zigzig (skew) diagram $D_A$ which
in French notation consists of rows of size $c_1, \ldots, c_{k+1}$, reading
from top to bottom, which overlap by one square. For example, if
$n=8$ and $A = \{2,4,7\}$, then $comp(A) = (2,2,3,1)$ and
$D(A)$ is the diagram pictured in Figure (\ref{fig:DA}).

\begin{df}
A bipartition (i.e., a pair of partitions) $\la=(\mu,\nu)$
where $\mu =(\mu_1 \geq \cdots \geq \mu_{k+1} \geq 0)$ and
$\nu = (\nu_1 \geq \cdots \geq \nu_{n-k+1} \geq 0)$
is called an $(n,k)$-bipartition if
\begin{enumerate}
\item $\lambda_k = \nu_{n-k+1} =0$ so that $\mu$ has at most $k-1$ parts and
$\nu$ has at most $n-k$ parts, .
\item for $i =1, \ldots, k-1$,
$\lambda_i - \lambda_{i+1} \in \{0,1\}$, and
\item  for $i =1, \ldots, n-k$,
$\nu_i - \nu_{i+1} \in \{0,1\}$.
\end{enumerate}

\noindent For a permutation $\pi\in S_n$ and a corresponding
$k$-descent basis element $a_\pi^{(k)} = \prod_{i=1}^{k-1}
x_{\pi(i)}^{d_i}\cdot \prod_{i=k+1}^{n} y_{\pi(i)}^{d_i}$, let
$$
\la(a_\pi^{(k)}) := ((d_1, d_2, \ldots, d_{k-1},0),(d_{n},d_{n-1},
\ldots, d_{k+1},0))
$$
be its {\em exponent bipartition}.

For an $(n,k)$-bipartition $\la=(\mu,\nu)$ let
$$
I_{\la}^{(k)\underline{\triangleleft}} :=
\sn_{\bbq}\{a^{(k)}_\pi+J_{(1^{k-1},n-k+1)}\,:\, \pi\in S_n,\
\la(a^{(k)}_\pi)\ \underline{\triangleleft}\ \la\,\},
$$
and
$$
I_{\la}^{(k)\triangleleft} :=
\sn_{\bbq}\{a^{(k)}_\pi+J_{(1^{k-1},n-k+1)}\,:\, \pi\in S_n,\
\la(a^{(k)}_\pi)\ {\triangleleft}\ \la\,\}
$$
be subspaces of the module ${\bf H}_{(1^{k-1},n-k+1)}'$, where
$\underline{\triangleleft}$ is the dominance order on bipartitions
(see Definition~\ref{order}.1), and let
$$
R_{\la}^{(k)}:= I_{\la}^{(k)\underline{\triangleleft}}/
I_{\la}^{(k)\triangleleft}.
$$
\end{df}

\begin{pro}
$I_{\la}^{(k)\underline{\triangleleft}}$,
$I_{\la}^{(k)\triangleleft}$ and thus $R_{\la}^{(k)}$ are
$S_n$-invariant.
\end{pro}


\begin{lem}\label{m-t.basis_R}
Let $\la=(\mu,\nu)$ be an $(n,k)$-bipartition. Then
\begin{equation}
\{a_{\pi}^{(k)}+I_{\la}^{(k)\triangleleft} \,:\, \Des(\pi) =
A_{\la}\}
\end{equation}
is a basis for $R_{\la}^{(k)}$, where
\begin{equation}\label{e.Rbasis3}
A_{\la} := \{ 1\leq i< n \, :\, \mu_i - \mu_{i+1} =1 {\rm\ or\ }
\nu_{n-i}-\nu_{n-i+1}=1\ \}. 
\end{equation}
\end{lem}

\begin{thm}\label{m-t.action}
The $S_n$-action on $R_{\la}^{(k)}$ is given by
$$
s_j(a_{\pi}^{(k)}) = \begin{cases} a_{s_j\pi}^{(k)}, &\mathrm{if\
} |\pi^{-1}(j+1) - \pi^{-1}(j)| > 1;\cr a_{\pi}^{(k)},
&\mathrm{if\ } \pi^{-1}(j+1) = \pi^{-1}(j)+1;\cr
-a_{\pi}^{(k)}-\sum_{\sigma\in A_j(\pi)} a_{\sigma}^{(k)},
  &\mathrm{if\ } \pi^{-1}(j+1) = \pi^{-1}(j)-1.\cr
\end{cases}
$$
Here $s_j = (j,j+1)$ $(1\leq j< n$) are the Coxeter generators of
$S_n$ and $\{a_{\pi}^{(k)}+I_{\la}^{(k)\triangleleft} \,:\,
\Des(\pi)=A_{\la}\}$
is the descent basis of $R_{\la}^{(k)}$. For the definition of
$A_j(\pi)$ see Theorem~\ref{t.action} below.
\end{thm}

This explicit description of the action is then used to prove the
following.

\begin{thm}
Let $\la=(\mu,\nu)$ be an $(n,k)$-bipartition. $R_{\la}^{(k)}$ is isomorphic,
as an $S_n$-module, to the Solomon descent representation determined by
the descent class $\{\pi\in S_n\,:\,\Des(\pi)=A_\la\}$.
\end{thm}



Let ${{\bf H}'}_{(1^{k-1},n-k+1)}^{(t_1,t_2)} $ be the $(t_1,t_2)$-th
homogeneous component of ${\bf H}_{(1^{k-1},n-k+1)}'$.

\begin{cor}\label{m-zigzags}
For every   $t_1,t_2 \geq 0$ and $1\leq k\leq n$, the
$(t_1,t_2)$-th homogeneous component of ${\bf H}_{(1^{k-1},n-k+1)}'$
decomposes into a direct sum of Solomon descent representations
as follows:
$$
{{\bf H}'}_{(1^{k-1},n-k+1)}^{(t_1,t_2)} \cong \bigoplus_{\la} R_{\la}^{(k)},
$$
where the sum is over all $(n,k)$-bipartitions $\la=(\mu,\nu)$.
such that
\begin{equation}\label{(n,k)weight}
 \sum\limits_{i\ge k{\rm\ and\ } \nu_i>\nu_{i+1}} (n-i) =t_1, \qquad
\sum\limits_{i<k{\rm\ and\ } \mu_i>\mu_{i+1}} i =t_2.
\end{equation}
\end{cor}

For example, suppose that $k=3$ and $n=4$. Then if
$\lambda =(\mu,\nu)$ is a $(4,3)$ partition, then
$\mu \in \{ (0,0,0),(1,0,0),(1,1,0),(2,1,0)\}$ and
$\nu \in \{ (0,0), (1,0)\}$. Table 1 below list
all the possible $(4,3)$-partitions. Then for
each $\lambda = (\mu,\nu)$, we  list  the corresponding
weight $q^{t_1}t^{t_2}$ given by (\ref{(n,k)weight}), the corresponding
descent set $A_\lambda$, and  the ribbon Schur function corresponding to
$A_\lambda$.

\begin{figure}[tbp]
\centering
$$
\epsfysize=4.0in \epsffile{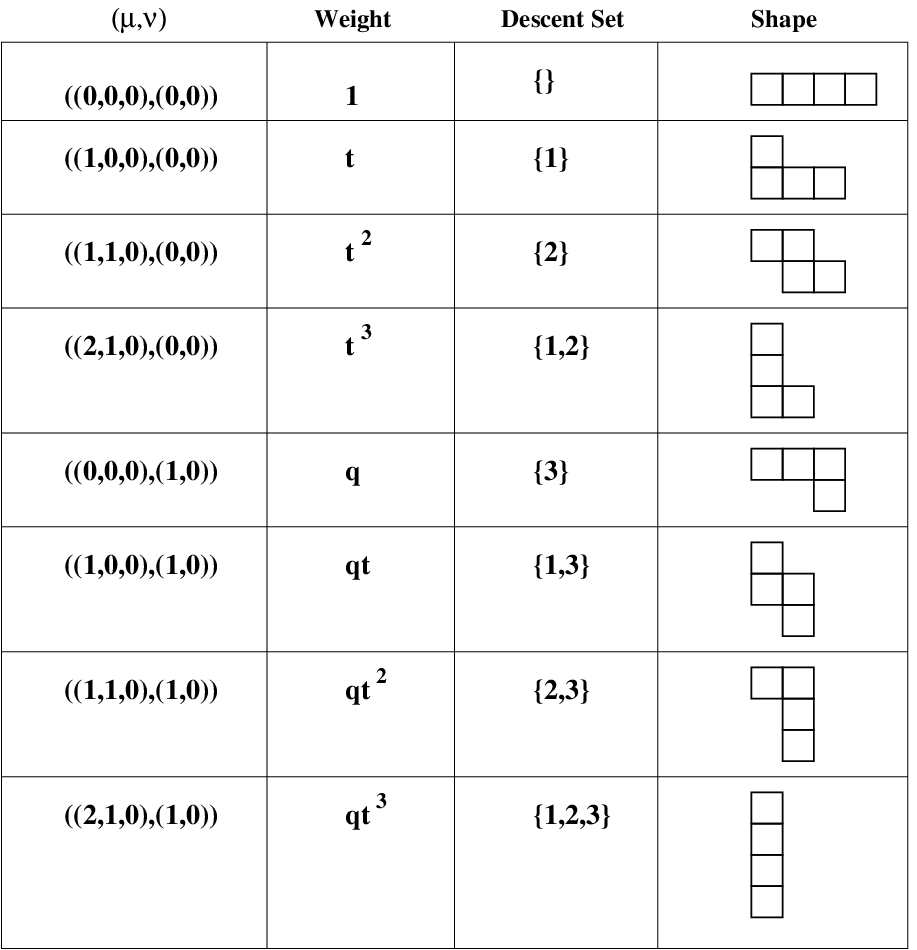}
$$
\caption{Table 1.}\label{fig:table}\end{figure}

\subsubsection{Decomposition into Irreducibles}\label{Intro-irreducibles}

A classical theorem of Lusztig and Stanley gives the multiplicity
of the irreducibles in the homogeneous component of the
coinvariant algebra of type $A$. Define $1\leq i<n$ to be a {\em
descent} in a standard Young tableau $T$ if $i+1$ lies strictly
above and weakly to the left of $i$ (in French notation). Denote
the set of all descents in $T$ by $\Des(T)$ and let the {\it major
index} of $T$ be
$$
\maj(T):=\sum\limits_{i\in \Des(T)} i .
$$

\begin{thm}\label{t.LS} \mbox{\rm\cite[Prop.~4.11]{St79}}
The multiplicity of the irreducible $S_n$-representation $S^\la$
in the $k$-th homogeneous component of the coinvariant algebra of
type $A$ is
$$
\#\{T\in SYT(\la)\,:\,\maj(T)=k\},
$$
where $SYT(\la)$ is the set of all standard Young tableaux of
shape $\la$.
\end{thm}

In 1994, Stembridge~\cite{Ste94} gave an explicit combinatorial
interpretation of the $(q,t)$-Kostka polynomials for hook shape.
Stembridge's result
implies the following extension of the Lusztig-Stanley theorem.

\begin{figure}[ht]
\begin{center}
\epsfig{file=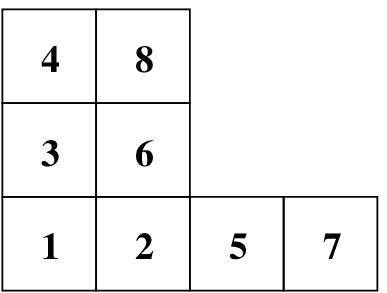} \caption{$\Des$, $\maj_{1,k}$, and
$\rmaj_{k,n}$ for a standard tableau.} \label{fig:Mac3}
\end{center}
\end{figure}

For a standard Young tableau $T$ define
$$
\maj_{i,j}(T) := \sum\limits_{r\in \Des(T) \atop i\leq r<j} r
$$
and
$$
\rmaj_{i,j}(T) := \sum_{r\in \Des(T) \atop i\leq r<j} (n-r).
$$
%
%
%
For example, for the column strict tableaux $T$ pictured in Figure
\ref{fig:Mac3} and $k=4$, $Des(T) = \{2,3,5,7\}$, $\maj_{1,4}(T) =
2 +3 =5$, and $\rmaj_{5,8} = (8-5)+(8-7) =4$.

We can restate Stembridge's results~\cite{Ste94} as follows.

\begin{thm} \label{t.stemb}
\begin{equation}
\tilde{K}_{\la,(1^{k-1},n-k+1)}(q,t) =
\sum_{T \in SYT(\la)} q^{\maj_{1,n-k+1}(T)} t^{\rmaj_{n-k+1,n}(T)}
\end{equation}
\end{thm}

Our next  result is an immediate consequence of Haiman's result
(\ref{char4}) and Theorem \ref{t.stemb}.

\begin{thm}\label{t.stemb2}
The multiplicity of the irreducible $S_n$-representation $S^\la$
in the $(h,h')$ level of  ${\bf H}_{(1^{k-1},n-k+1)}$ (bi-graded by
total degrees in the $x$'s and $y$'s) is
$$
\chi_{\la}^{(h,h')} = \#\{T\in SYT(\la)\,:\,
\maj_{1,n-k+1}(T)=h,\,\rmaj_{n-k+1,n}(T)=h'\},
$$
where $SYT(\la)$ is the set of all standard Young tableaux of
shape $\la$.
\end{thm}

Haglund~\cite{Hag04} gave another proof of Theorem~\ref{t.stemb2}
that used his conjectured combinatorial definition of
$\tilde{H}_\mu(\ox;q,t)$. Haglund's conjecture has recently been
proved by Haglund, Haiman and Loehr \cite{HHL05a,HHL05b}.

We give two proofs of this decomposition rule. The first one,
given in Section~\ref{s.des-rep}, follows from the decomposition
into descent representations described in
Corollary~\ref{m-zigzags} above.

The second proof, 
given in Section~\ref{s.decomposition}, is more 
``combinatorial". It  uses the mechanism of~\cite{HHL05b} but does
not rely on Haglund's combinatorial interpretation of
$\tilde{H}_{(1^{k-1},n-k+1)}(\ox;q,t)$.

\section{The Garsia-Haiman Module ${\bf H}_\mu$}


In this section, we shall provide some background on the
Garsia-Haiman module, also known as the space of
orbit harmonics. 
Proofs of the results in this section can be found in \cite{GH1}
and \cite{GHOH}.

Let $m:=2n$ and $(z_1,\ldots,z_m) := (x_1,\ldots,x_n,y_1,\ldots,y_n)$.
Given positive weights $w_1, \ldots, w_m$,
define a grading of ${\bf R}=\bbq[z_1,\dots,z_m]$ by setting,
for any monomial $\oz^p = z_1^{p_1} \cdots z_m^{p_m}$,
\begin{equation}\label{degree}
d_w(\oz^p) := \sum_{i=1}^m w_i p_i.
\end{equation}
For any polynomial $g = \sum_p c_p \oz^p$ and any $k$, set
\begin{equation}\label{component}
\pi^w_k g = \sum_{d_w(\oz^p) =k} c_p \oz^p
\end{equation}
Call $\pi^w_k g$ the {\it $w$-homogeneous component of $w$-degree
$k$}, and let $\htop_w(g)$ denote the $w$-homogeneous component of
highest $w$-degree in $g$.  Let $V$
be a subspace of ${\bf R}$. We say that $V$ is {\it
$w$-homogeneous} if $\pi_k^w V \subseteq V$ for all $k$. Let
$V^{\perp}$ is the orthogonal complement of $V$ with respect to
the inner product defined at the beginning of
Subsection~\ref{s.main}. If V is $w$-homogeneous, then $V^{\perp}$
is also $w$-homogeneous.

If $J$ is an ideal of ${\bf R}$ then
\begin{equation}
J^{\perp} = \{g \in \bbq[\oz]:
f(\frac{\partial}{\partial z_1},\ldots,\frac{\partial}{\partial z_m})g(\oz) = 0
\ (\forall f \in J)\}.
\end{equation}
In particular, $J^{\perp}$ is closed under differentiation. We set
${\bf R}_J = {\bf R}/J$ and define the {\em associated $w$-graded
ideal of $J$} to be the ideal
\begin{equation}\label{gr}
gr_w J := \langle \htop_w(g):g \in J\rangle.
\end{equation}
The quotient ring ${\bf R}/gr_w J$ is referred to as the {\em
$w$-graded version of ${\bf R}_J$} and is denoted by $gr_w {\bf
R}_J$.

Given a point $\rho =(\rho_1, \ldots, \rho_m)$, denote its
$S_n$-orbit by
\begin{equation}\label{orbit}
[\rho] = [\rho]_{S_n} = \{\sigma(\rho) : \sigma \in S_n\}.
\end{equation}
We also set
\begin{equation}
{\bf R}_{[\rho]}: = {\bf R}/J_{[\rho]}
\end{equation}
where
\begin{equation}
J_{[\rho]} = \{g \in {\bf R}: g(\bar z) =0 \ (\forall \bar z \in
[\rho])\}
\end{equation}
Considered as an algebraic variety, ${\bf R}_{[\rho]}$ is the
coordinate ring of $[\rho]$. Since the ideal $J_{[\rho]}$ is
$S_n$-invariant, $S_n$ also acts on ${\bf R}_{[\rho]}$. In fact,
it is easy to see that the corresponding representation is
equivalent to the action of $S_n$ on the left cosets of the
stabilizer of $\rho$. Thus if the stabilizer of $\rho$ is trivial,
i.e. if $\rho$ is a regular point, then ${\bf R}_{[\rho]}$ is a
version of the left regular representation of $S_n$. If each
$\sigma \in S_n$ preserves $w$-degrees, then we can associate two
further $S_n$-modules with $\rho$; namely, $gr_w J_{[\rho]}$ and
its orthogonal complement
\begin{equation}\label{harmonics}
{\bf H}_{[\rho]}= (gr_w J_{[\rho]})^\perp.
\end{equation}
Garsia and Haiman refer to ${\bf H}_{[\rho]}$ as the space of {\em
harmonics of the orbit $[\rho]$}. Clearly if $f(\oz)$ is an
$S_n$-invariant polynomial, then $f(\oz) -f(\rho) \in J_{[\rho]}$.
In addition, if $f(\oz)$ is $w$-homogeneous, then $f(\oz) \in gr_w
J_{[\rho]}$. This implies that any element $g \in {\bf
H}_{[\rho]}$ must satisfy the differential equation
$$
f(\frac{\partial}{\partial z_1}, \ldots ,\frac{\partial}{\partial z_m})
g(z_1, \ldots, z_m) = 0.
$$
It is not difficult to show that ${\bf H}_{[\rho]}$ and $gr_w {\bf
R}_{[\rho]}$ are equivalent $w$-graded modules and that these two
spaces realize a graded version of the regular representations of
$S_n$.

Given a partition $\mu = (\mu_1 \geq \cdots \geq \mu_k > 0)$ with
$k$ parts, let $h =\mu_k$ be the number of parts of the conjugate
partition $\mu^\prime$. Let $\alpha_1, \ldots, \alpha_k$ and
$\beta_1, \ldots, \beta_h$ be distinct rational numbers.
Alternatively, we can think of the $\alpha$'s and $\beta$'s as
indeterminates. An {\it injective tableau of shape $\mu$} is a
labelling of the cells of $\mu$ by the numbers $\{1,2,\ldots,n\}$
so that no two cells are labelled by the same number. The
collection of all such tableaux will be denoted by ${\mathcal
IT}_\mu$. For each $T \in {\mathcal IT}_\mu$ and $\ell \in
\{1,2,\ldots,n\}$, let $s_T(\ell) = (i_T(\ell),j_T(\ell))$ denote
the cell of $\mu$ which contains the number $\ell$. 
We then
construct a point $\rho(T) = (a(T),b(T))$ in $2n$-dimensional
space by setting
\begin{equation}\label{abT}
a_\ell(T) := \alpha_{i_T(\ell)} \ \ \mbox{and} \ \ b_\ell(T) =
\beta_{j_T(\ell)} \qquad (1\le \ell\le n).
\end{equation}

\begin{figure}[ht]
\begin{center}
\epsfig{file=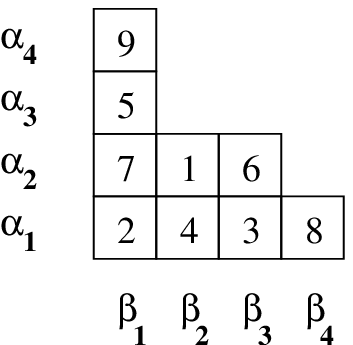} \caption{$\rho(T)$.} \label{fig:abT}
\end{center}
\end{figure}

For example, for the injective tableau $T$ pictured in Figure
\ref{fig:abT},
$$\rho(T) =
(\alpha_2,\alpha_1,\alpha_1,\alpha_1,\alpha_3,\alpha_2,\alpha_2,\alpha_1,\alpha_4;
\beta_2,\beta_1,\beta_3,\beta_2,\beta_1,\beta_3,\beta_1,\beta_4,\beta_1)$$

The collection $\{\rho(T): T \in {\mathcal IT}_\mu\}$ consists of
$n!$ points. In fact, it is simply the $S_n$-orbit under the
diagonal action of any one of its elements, where the {\em
diagonal action} of $S_n$ on $2n$-dimensional space is defined by
$$
\sigma(x_1, \ldots, x_n;y_1, \ldots, y_n) := (x_{\sigma_1},
\ldots, x_{\sigma_n};y_{\sigma_1}, \ldots, y_{\sigma_n})
$$
for any $\sigma \in S_n$. For the rest of this paper we let
$$
{\bf R} = \bbq[\ox,\oy] = \bbq [x_1, \ldots, x_n;y_1, \ldots, y_n]
$$
and use the grading $w_i=1$ for all $i=1, \ldots, 2n$. We shall
let ${\bf R}_{[\rho_\mu]}$, $gr {\bf R}_{[\rho_\mu]}$, and ${\bf
H}_{[\rho_\mu]}$ denote the corresponding coordinate ring, its
graded version, and the space of harmonics for the $S_n$-orbit
$\{\rho(T): T \in {\mathcal IT}_\mu\}$.

Garsia and Haiman \cite{GH1} proved the following.

\begin{pro}
If $(i,j)$ is a cell outside of $\mu$, then for any $s \in \{1,
\ldots, n\}$, the monomial $x_s^{i-1}y_s^{j-1}$ belongs to the
ideal $gr J_{[\rho_\mu]}$. In particular, if a monomial
$x^py^q = x_1^{p_1} \cdots x_n^{p_n}y_1^{q_1}\cdots
y_n^{q_n}\not\in J_{[\rho_\mu]}$, then all the pairs $(p_s,q_s)$ must be
cells
of $\mu$.
\end{pro}

\begin{thm}\label{gh}
Let ${\bf H}_\mu$ be the vector space of polynomials spanned by
all the partial derivatives of $\Delta_\mu(x_1, \ldots,x_n;y_1,
\ldots, y_n)$ (see (\ref{Dmu}) above). Then
\begin{enumerate}
\item ${\bf H}_\mu \subseteq {\bf H}_{[\rho_\mu]}$.

\item If ${\rm dim}({\bf H}_\mu) = n!$, then
${\bf H}_\mu = {\bf H}_{[\rho_\mu]}$ and
$J_\mu = {\bf H}_\mu^\perp = gr J_{[\rho_\mu]}$.
\end{enumerate}

\end{thm}

\bigskip

In 2001, Haiman solved the $n!$ conjecture and proved that the
dimension of ${\bf H}_\mu$ is indeed $n!$~\cite{Hai01}. Thus,
Theorem~\ref{gh} implies that for every partition $\mu$,
$J_\mu = {\bf H}_\mu^\perp$ is the ideal in ${\bf R}$ generated by the set
$$
gr\,\{f\in {\bf R}\, :\, f(\rho(T))=0\ \ (\forall T\in {\mathcal
IT}_\mu) \}.
$$

We deduce

\begin{cor}\label{s0}
The following polynomials belong to the ideal $J_{(1^{k-1},n-k+1)}$:
\begin{itemize}
\item[(i)]
$\Lambda(x_1,\dots,x_n)^+$ and
$\Lambda(y_1,\dots,y_n)^+$ (the symmetric functions in $x$ and $y$
without a constant term),
\item[(ii)]
the monomials $x_{i_1}\cdots x_{i_k}$ ($i_1<\cdots<i_k$) and
$y_{i_1}\cdots y_{i_{n-k+1}}$ ($i_1<\cdots<i_{n-k+1}$), and
\item[(iii)]
the monomials $x_iy_i$ ($1\le i\le n$).
\end{itemize}
\end{cor}

\noindent{\bf Proof.}
\begin{itemize}
\item[(i)]
The multiset of row-coordinates is the same for all injective tableaux of
a fixed shape. Thus, for any $f\in \Lambda(x_1,\dots,x_n)^+$,
$g(\bar z) := f(x_1,\dots,x_n)-f(\alpha_1,\dots,\alpha_1,\alpha_2,\dots,\alpha_k)$
is zero for all substitutions $\bar z = \rho(T)$, where $T$ is an injective
tableau of shape $(n-k+1,1^{k-1})$. Thus $f=\htop(g)\in J_{(n-k+1,1^{k-1})}$.
Similarly, for any $f\in \Lambda(y_1,\dots,y_n)^+$,
$g(\bar z) := f(y_1,\dots,y_n)-f(\beta_1,\dots,\beta_1,\beta_2,\dots,\beta_{n-k+1})$
is zero for all substitutions $\bar z = \rho(T)$. Thus $f=\htop(g)\in J_{(n-k+1,1^{k-1})}$.
\item[(ii)]
In the shape $(1^{k-1},n-k+1)$ there are only $k-1$ cells
outside of row $1$. Therefore
$g(\bar z) := (x_{i_1}-\alpha_1)\cdots(x_{i_k}-\alpha_1)$ is zero
for all $1\le i_1<\cdots <i_k\le n$ and all $\bar z = \rho(T)$.
Thus $x_{i_1}\cdots x_{i_k} = \htop(g)\in J_{(n-k+1,1^{k-1})}$.
Similarly,
$g(\bar z) := (y_{i_1}-\beta_1)\cdots(y_{i_{n-k+1}}-\beta_1)$ is zero
for all $1\le i_1<\cdots <i_{n-k+1}\le n$ and all $\bar z = \rho(T)$.
Thus $y_{i_1}\cdots y_{i_{n-k+1}} = \htop(g)\in J_{(n-k+1,1^{k-1})}$.
\item[(iii)]
Each cell of a hook shape is either in row $1$ or in column $1$. Thus
$g(\bar z) := (x_i-\alpha_1)(y_i-\beta_1)$ is zero for all $1\le i\le n$ and all
$\bar z = \rho(T)$, so that $x_i y_i = \htop(g)\in J_{(n-k+1,1^{k-1})}$.
\end{itemize}

\qed

\medskip

In the sequel we will show that the polynomials in
Corollary~\ref{s0} actually generate the ideal $J_{(1^{k-1},n-k+1)}$.

%

\section{Generalized Kicking-Filtration Process}\label{s.kick}

In this section we prove Theorem~\ref{mk}, which  implies
Corollary~\ref{m11} as a special
case and Theorem~\ref{t.m2}. The idea is to generalize the kicking process 
for obtaining a basis. The kicking process was used in an early
paper of Garsia and Haiman~\cite{GH1} to prove the $n!$-conjecture
for hooks. We combine ingredients of this process with a
filtration.

\subsection{Proof of Theorem~\ref{mk}}\label{pf-mk}


%
%
%
%
%
For every triple $(A,c,\bar A)$, where $[n]=A\cup \{c\}\cup \bar
A$ and $|A|=k-1, |\bar A|=n-k$, define an $(A,c,\bar
A)$-permutation $\pi_{(A,c,\bar A)}\in S_n$, in which the letters
of $A$ appear in decreasing order, then $c$, and then the
remaining letters in increasing order. For example, if $n=9, k=4,
c=5, A=\{1,6,7\}$ then
 $
 \pi_{(A,c,\bar A)}=761523489$.


For given $n$ and $k$, order the $N:=n{n-1\choose k-1}$ distinct
$(A,c,\bar A)$-permutations in reverse lexicographic order (as words):
$\pi_1,\dots,\pi_N$. If $\pi_i$ corresponds to $(A,c,\bar A)$, let
$m_i:=m_{(A,c,\bar A)}$ (see Subsection~\ref{s.k-AH} above).

For example, for $n=4$ and $k=3$, the list of $(A,c,\bar A)$-permutations is
$$
 \pi_{(\{34\},2,\{1\})}=4321,\
 \pi_{(\{34\},1,\{2\})}=4312,\
 \pi_{(\{24\},3,\{1\})}=4231,\ldots,
$$
$$
 \pi_{(\{13\},2,\{4\})}=3124,\
  \pi_{(\{12\},4,\{3\})}=2143,\
 \pi_{(\{12\},3,\{4\})}=2134.
$$
and the order is
$$
4321<_L  4312<_L  4231 <_L  4213 <_L 4132  <_L 4123 <_L 3241<_L 3214<_L
3142<_L 3124 <_L 2143 <_L 2134.
$$
Thus the permutations are indexed by
$\pi_1=4321,\ \pi_2=4312,\ \pi_3=4231,\ \dots, \pi_{11}=2143,\ \pi_N=\pi_{12}=2134$
and the corresponding monomials are
$m_1= x_4 x_3 y_1,\ m_2=x_4x_3,\ m_3= x_4 y_1, \ldots, \ m_{11}= y_3,\ m_N=m_{12}=1$.



\bigskip

Let
$$
I_0:=J_{(1^{k-1},n-k+1)}
$$
and define recursively
$$
I_t:= I_{t-1}+ m_t \bbq[\bar x,\bar y] \qquad (1\le t \le N).
$$
Clearly,
$$
I_0\subseteq I_1\subseteq I_2\subseteq \cdots \subseteq
I_N=\bbq[\bar x,\bar y].
$$
The last equality follows from the fact that $m_N=1$.

\begin{obs}
$$
{\bf H}_{(1^{k-1},n-k+1)}'=\bbq[\bar x,\bar y]/I_0\cong
\bigoplus_{t=1}^N (I_t/I_{t-1})
$$
as vector spaces. In particular, a sequence of bases for the
quotients $I_t/I_{t-1}$, $1\le t\le N$, will give a basis for
${\bf H}_{(1^{k-1},n-k+1)}'$.
\end{obs}

It remains to prove that the set $m_t\cdot B_A\cdot C_{\bar A}$,
where $B_A,\ C_{\bar A}$ are bases of coinvariant algebras in
$\bar x_A$ and $\bar y_{\bar A}$ respectively, is a basis for
$I_t/I_{t-1}$.



\smallskip

Consider the natural projection
$$
f_t: m_t \bbq[\bar x,\bar y]\longrightarrow I_t/I_{t-1} .
$$
Clearly, $f_t$ is a surjective map. 

\begin{lem}\label{t.kick-claim}
$$
m_t\cdot \left(\sum\limits_{i\not\in A} \langle x_i \rangle +
\sum\limits_{j\not\in \bar A} \langle y_j \rangle + \langle
\Lambda[\bar x]^+ \rangle + \langle \Lambda[\bar y]^+
\rangle\right) \ \subseteq \ I_{t-1}\cap m_t \bbq[\bar x, \bar y]
= \kr (f_t),
$$
so that $f_t$ is well defined on the quotient
$$
m_t \bbq[\bar x,\bar y]/ \left[m_t\cdot\left(\sum\limits_{i\not\in
A} \langle x_i \rangle + \sum\limits_{j\not\in \bar A} \langle y_j
\rangle + \langle \Lambda[\bar x]^+ \rangle + \langle \Lambda[\bar
y]^+ \rangle\right)\right]\cong $$ $$ m_t \cdot \bbq[\bar
x_A]/\langle\Lambda[\bar x_A]^+\rangle \cdot \bbq[\bar y_{\bar
A}]/\langle\Lambda[\bar y_{\bar A}]^+\rangle).
$$
\end{lem}

\noindent{\bf Proof of Lemma~\ref{t.kick-claim}.} First, let
$i\not\in A$. We shall show that $ m_t x_i \in I_{t-1}$.
Consider the four possible cases:
\begin{itemize}
\item[Case X1.] $a>c$ for all $a\in A$.\\ Then $m_t$ has all the
$x_a$ ($a\in A$) as factors, so that $m_t x_i$ has exactly $k$
distinct $x$-factors and therefore belongs to $I_0\subseteq
I_{t-1}$ (by Corollary~\ref{s0}(ii)).

\item[Case X2.] $i<c$ and there exists $a\in A$ such that $a<c$.\\
$i\not\in A$ and $i\ne c$, hence $i\in \bar A$. As $i<c$, $y_i$ is
a factor of $m_t$. Hence $x_i y_i$ is a factor of $m_t x_i$
 and therefore belongs to $I_0\subseteq I_{t-1}$ (by Corollary~\ref{s0}(iii)).

\item[Case X3.] $i=c$ and there exists $a\in A$ such that $a<c$.\\
Let $a_0:= \max\{a\in A\, :\, a<c\}$. Let $t_{(c,a_0)}$ be the
transposition interchanging $c$ with $a_0$ and let $A':=
(A\setminus\{a_0\})\cup \{c\}$. By the definition of $a_0$,
$m_{t'}=m_{(A',a_0,\bar {A'})}$
 divides  $m_t x_i$. To verify this, notice
that for every $j$, such that $x_j$ divides $m_{t'}$, $j\in A'$
and $j>a_0$. Hence $j\ne a_0$ and thus $j\in A$. Also, by the
definition of $a_0$, for every $j\in A$, $j>a_0 \Rightarrow j> c$.
It follows that $x_j$ divides $m_t$.
 Finally, for every $j$, such $y_j$ divides
$m_{t'}$, $j\in \bar A= \bar {A'}$ and $j<a_0$. Hence $j<c$ and
$y_j$ divides $m_t$.

On the other hand, $\pi_{(A',a_0,\bar {A'})}=
t_{(c,a_0)}\pi_{(A,c,\bar A)} <_L \pi_{(A,c,\bar A)}$, since
$a_0<c$. Hence $m_{t'}=m_{(A',a_0,\bar {A'})}\in I_{t-1}$ and
$m_t x_i\in \langle m_{t'}\rangle\subseteq I_{t-1}$.

\item[Case X4.] $i>c$ and there exists $a\in A$ such that $a<c$.\\
Let $a_0:= \max\{a\in A\, :\, a<c\}$ as above. Let $t_{(a_0,i,c)}$
be the 3-cycle mapping $a_0$ to $i$, $i$ to $c$ and $c$ to $a_0$.
Let $A':= (A\setminus\{a_0\})\cup \{i\}$ and
$m_{t'}=m_{(A',a_0,\bar{A'})}$. Then $\pi_{(A',a_0,\bar
{A'})}\le_L t_{(a_0,i,c)}\pi_{(A,c,\bar A)} <_L \pi_{(A,c,\bar
A)}$, since $a_0<c<i$. Thus $t'<t$.

On the other hand, $m_{t'}=m_{(A',a_0,\bar {A'})}$
 divides  $m_t x_i$. This follows from the  implications below:
 $$
i'\in A' \hbox{ and } i'>a_0 \Longrightarrow i'=i \hbox{ or } i'>c
.
$$
Hence, every $x_{i'}$ which divides $m_{t'}$ divides also $m_t
x_i$. Also
$$
j'\in \bar {A'} \hbox{ and } j'<a_0 \Longrightarrow j'\in \bar A
\hbox{ and } j'< c ,
$$
Hence, every $y_{j'}$ which divides $m_{t'}$ divides also $m_t
x_i$.

We conclude that $m_{t'}=m_{(A',a_0,\bar A)}\in I_{t-1}$ and
$\langle m_t x_i\rangle\subseteq \langle m_{(A',a_0,\bar
A)}\rangle \subseteq I_{t-1}$.

\end{itemize}

Similarly, by considering four analogous cases, one can show that
if $j\not\in \bar A$ then $ m_t y_j \in I_{t-1}$.

\medskip

Finally, observe that, by Corollary~\ref{s0}(i), $m_t \langle
\Lambda[\bar x]^+\rangle \subseteq I_0\subseteq I_{t-1}$ and $m_t
\langle \Lambda[\bar y]^+\rangle \subseteq I_0\subseteq I_{t-1}$.

\medskip

This completes the proof of Lemma~\ref{t.kick-claim}.
\qed

\bigskip

Now, recall that $B_A$ and $C_{\bar A}$ are bases for the
coinvariant algebras $\bbq[\bar x_A]/\langle \Lambda[\bar
x_A]^+\rangle$ and $\bbq[\bar y_{\bar A}]/\langle\Lambda[\bar
y_{\bar A}]^+\rangle$ respectively,
 so
that $m_t\cdot B_A\cdot C_{\bar A}$ 
spans $m_t\cdot \bbq[\bar x_A]/\langle \Lambda[\bar x_A]^+\rangle
\cdot \bbq[\bar y_{\bar A}]/\langle\Lambda[\bar y_{\bar
A}]^+\rangle$.

In order to prove Theorem~\ref{mk}, it remains to show that for
every $1\le t\le N$,
\begin{itemize}
\item[(a)] $m_t\cdot B_A\cdot C_{\bar A}$ is a basis for $m_t\cdot
\bbq[\bar x_A]/\langle \Lambda[\bar x_A]^+\rangle \cdot \bbq[\bar
y_{\bar A}]/\langle\Lambda[\bar y_{\bar A}]^+\rangle$, and
\item[(b)] $f_t$ is one-to-one.
\end{itemize}
Indeed, by Lemma~\ref{t.kick-claim},
$$
\dm (I_t/I_{t-1})\le \dm (m_t \cdot
\bbq[\bar x_A]/\langle\Lambda[\bar x_A]^+\rangle \cdot
\bbq[\bar y_{\bar A}]/\langle\Lambda[\bar y_{\bar A}]^+\rangle) \le
$$
\begin{equation}\label{e.dim_quotient}
\dm (\bbq[\bar x_A]/\langle\Lambda[\bar x_A]^+\rangle
     \bbq[\bar y_{\bar A}]/\langle\Lambda[\bar y_{\bar A}]^+\rangle)
= (k-1)! \cdot (n-k)!
\end{equation}
The last equality follows from the classical result that the
coinvariant algebra of $S_n$ carries the regular
$S_n$-representation~\cite[\S 3]{Hum}.

If there exists $1\le t\le N$ such that either (a) or (b) does not hold,
then there exists $t$ for which a sharp inequality holds
in~(\ref{e.dim_quotient}). Then
$$
\dm {\bf H}_{(1^{k-1},n-k+1)}'=\dm \left(\bbq[\bar x,\bar
y]/I_0\right)= \sum\limits_{t=1}^N \dm (I_t/I_{t-1})< N (k-1)!
(n-k)! = n  {n-1\choose k-1} (k-1)! (n-k)! = n! ,
$$
contradicting the $n!$ theorem. This completes the proof of
Theorem~\ref{mk}.

\qed

\subsection{Applications}

\noindent{\bf Proof of Corollary~\ref{m11}.} By choosing $B_A$ and
$C_{\bar A}$ in Theorem~\ref{mk} as Artin bases of the
corresponding coinvariant algebras we get the $k$-th Artin basis
($b^{(k)}_\pi$).

By choosing $B_A$ and $C_{\bar A}$ as descent bases  we get the
$k$-th descent basis ($a^{(k)}_\pi$).

Finally, by choosing $B_A$ as a descent basis and $C_{\bar A}$ as
an Artin basis we get the $k$-th Haglund basis ($c^{(k)}_\pi$).

\qed

\bigskip

\noindent{\bf Proof of Theorem~\ref{t.m2}.} Let $J'$ be the ideal
of $\bbq[\bar x,\bar y]$ generated by
 $\Lambda(x_1,\dots,x_n)^+$,
$\Lambda(y_1,\dots,y_n)^+$ and the monomials $x_{i_1}\cdots
x_{i_k}$ ($i_1<\cdots<i_k$), $y_{i_1}\cdots y_{i_{n-k+1}}$
($i_1<\cdots<i_{n-k+1}$), and $x_iy_i$ ($1\le i\le n$). By
Corollary~\ref{s0}, $J'\subseteq J_{(1^{k-1},n-k+1)}$.

In the above proof of Theorem~\ref{mk}, if one replaces
$I_0:=J_{(1^{k-1},n-k+1)}$ by $I_0:=J'$, the same conclusions hold:
$$
\dm \left(\bbq[\bar x,\bar y]/J'\right)= \sum\limits_{t=1}^N \dm
(I_t/I_{t-1}) \le N (k-1)!(n-k)!= n!
$$
On the other hand, by~\cite{GH1},
$$
\dm \left(\bbq[\bar x,\bar y]/J'\right)\ge \dm \left(\bbq[\bar
x,\bar y]/J_{(1^{k-1},n-k+1)}\right)= \dm {\bf H}_{(1^{k-1},n-k+1)}'= n!
$$
Thus equality holds everywhere, and $J'= J_{(1^{k-1},n-k+1)}$. \qed

%


\bigskip

\section{A $k$-th Analogue of the Polynomial Ring}

The current section provides an appropriate setting for an
extension of the straightening algorithm from the coinvariant
algebra~\cite{Allen1, TAMS} to the Garsia-Haiman hook modules. The
algorithm will be given in Section~\ref{s.straight} and will be
used later to describe the $S_n$-action on ${\bf
H}_{(1^{k-1},n-k+1)}'$ and resulting decomposition rules (see
Section~\ref{s.des-rep}).

\subsection{$\DP_n^{(k)}$ and its Monomial Basis}

\begin{df}\label{3.1}
For every $1\le k\le n$ let ${\cal I}_k$ be the ideal in
$\bbq[x_1,\dots,x_n,y_1,\dots,y_n]$ generated by
\begin{itemize}
\item[(i)]
the monomials $x_{i_1}\cdots x_{i_k}$ $(1\le i_1<\cdots<i_k\le n)$,
\item[(ii)]
the monomials $y_{i_1}\cdots y_{i_{n-k+1}}$ $(1\le i_1<\cdots<i_{n-k+1}\le n)$, and
\item[(iii)]
the monomials $x_i y_i$ $(1\le i\le n)$.
\end{itemize}
Denote
$$
\DP_n^{(k)}:=\bbq[x_1,\dots,x_n,y_1,\dots,y_n]/{\cal I}_k.
$$
\end{df}

\begin{cla}\label{32}
The ideal ${\cal I}_k$ is contained in the ideal
$J_{(1^{k-1},n-k+1)}$.
\end{cla}

\noindent{\bf Proof.} Immediate from Definition~\ref{3.1} and
Corollary~\ref{s0}. \qed

%

\begin{df}
For a monomial $m=\prod\limits_{i=1}^n
x_i^{p_i}\prod\limits_{j=1}^n y_j^{q_j}\in
\bbq[x_1,\dots,x_n,y_1,\dots,y_n]$
define the {\em $x$-support} and the {\em $y$-support}
$$
\supp_x(m):= \{i\,:\,p_i > 0\}, \qquad \supp_y(m):= \{j\,:\,q_j >
0\}.
$$
\end{df}

\begin{df}
Let $M_n^{(k)}$ be the set of all monomials in
$\bbq[x_1,\dots,x_n,y_1,\dots,y_n]$ with
\begin{itemize}
\item[(i)]
$|\supp_x(m)| \le k-1$
\item[(ii)]
$|\supp_y(m)| \le n-k$, and
\item[(iii)]
$\supp_x(m)\cap \supp_y(m) = \emptyset$.
\end{itemize}
\end{df}

\begin{obs}
$\{m+{\cal I}_k\,:\,m\in M_n^{(k)}\}$ is a basis for
$\DP_n^{(k)}$.
\end{obs}

\begin{exa}
Let $k=n$. Then
$$
\DP_n^{(n)} \cong \bbq[x_1,\dots,x_n]/\langle x_1\cdots
x_n\rangle,
$$
and $M_n^{(n)}$ consists of all the monomials in $\bbq[x_1,\dots,x_n]$
which do not involve all the $x$ variables, i.e., are not divisible by
$x_1\cdots x_n$.
Similarly,
$$
\DP_n^{(1)} \cong \bbq[y_1,\dots,y_n]/\langle y_1\cdots
y_n\rangle,
$$
and $M_n^{(1)}$ consists of all the monomials in $\bbq[y_1,\dots,y_n]$
which do not involve all the $y$ variables.
\end{exa}

\subsection{A Bijection}

Next, define a map $\psi^{(k)}: M_n^{(k)}\to M_n^{(n)}$.

Every monomial $m\in M_n^{(k)}$ has the form $m =
x_{i_1}^{p_{i_1}}\cdots x_{i_{k-1}}^{p_{i_{k-1}}}\cdot
y_{j_1}^{q_{j_1}}\cdots y_{j_{n-k}}^{q_{j_{n-k}}}$ (with disjoint
supports of $x$'s and $y$'s). Let $u:=\max_j q_j$ and define
$$
\psi^{(k)}(m) := x_{i_1}^{p_{i_1}}\cdots
x_{i_{k-1}}^{p_{i_{k-1}}}\cdot x_{j_1}^{-q_{j_1}}\cdots
x_{j_{n-k}}^{-q_{j_{n-k}}}\cdot (x_1\cdots x_n)^u.
$$
Notice that if $u=0$, then $\psi^{(k)}(m) = m\in M_n^{(k)}\cap \bbq[x_1,\ldots,x_n]$ has
$|\supp_x(m)|\le k-1\le n-1$ so that $\psi^{(k)}(m)\in M_n^{(n)}$,
and, if $u>0$ is attained at $y_{j_0}$, then the exponent of $x_{j_0}$ in
$\psi^{(k)}(m)$ is zero so that again $\psi^{(k)}(m)\in M_n^{(n)}$.

\begin{exa}
Let $n=8$, $k=5$, and $m = x_3^6 x_4^6 x_7^5 y_6^2 y_2^3\in M^{(5)}_8$.
Then $u=3$, $j_0=2$, and
$$
\psi^{(k)}(m) = x_3^6 x_4^6 x_7^5 x_6^{-2} x_2^{-3} \cdot (x_1\cdots x_8)^3
= x_3^9 x_4^9 x_7^8 x_1^3 x_5^3 x_8^3 x_6^1\in M_8^{(8)}.
$$
\end{exa}

We claim that the map $\psi^{(k)}: M_n^{(k)}\to M_n^{(n)}$ is a bijection.
This will be proved by defining an inverse map $\phi^{(k)}: M_n^{(n)}\to M_n^{(k)}$.

Recall from \cite{TAMS} that the {\it index permutation} of a monomial
$m = \prod_{i=1}^n x_i^{p_i}\in \bbq[x_1,\dots,x_n]$ is the unique
permutation $\pi=\pi(m) \in S_n$ such that
$$
p_{\pi(i)}\ge p_{\pi(i+1)} \qquad (1\le i < n)  \leqno(i)
$$
and
$$
p_{\pi(i)}=p_{\pi(i+1)} \Longrightarrow \pi(i)<\pi(i+1).
\leqno(ii)
$$
In other words, $\pi$ reorders the variables $x_i$ by (weakly) decreasing
exponents, where the variables with a given exponent are ordered by
increasing indices.
Also, let $\la = \la(m) = (\la_1,\ldots,\la_n)$ be the corresponding
{\em exponent partition}, where $\la_i$ is the exponent of $x_{\pi(i)}$ in $m$.
Finally, let $\la'=(\la'_1,\dots,\la'_\ell)$ be the partition
conjugate to $\la$.

Every monomial $m = \prod_{i=1}^{n} x_{\pi(i)}^{\la_i}\in M_n^{(n)}$
can thus be written in the form
$$
m = \prod_{t=1}^\ell (x_{\pi(1)}\cdots x_{\pi(\la'_t)}).
$$
Note that, since $m$ is not divisible by $x_1 \cdots x_n$, $\la_n = 0$
and therefore $\la'_t < n$ $(\forall t)$.
Define
$$
\phi^{(k)}(m) :=
\prod_{\la'_t\le k-1} (x_{\pi(1)}\cdots x_{\pi(\la'_t)}) \cdot
\prod_{\la'_t\ge k} (y_{\pi(\la'_t+1)}\cdots y_{\pi(n)}).
$$

\begin{pro}\label{t.inverse_maps}
$\psi^{(k)}: M_n^{(k)}\to M_n^{(n)}$ and $\phi^{(k)}: M_n^{(n)}\to M_n^{(k)}$
are inverse maps.
\end{pro}
\noindent
{\bf Proof.\ }
Take a monomial $m = \prod_{i=1}^{n} x_{\pi(i)}^{\la_i}\in M_n^{(n)}$
as above. Then
\begin{eqnarray*}
\supp_x(\phi^{(k)}(m)) &\subseteq& \{\pi(1),\ldots,\pi(k-1)\},\cr
\supp_y(\phi^{(k)}(m)) &\subseteq& \{\pi(k+1),\ldots,\pi(n)\},
\end{eqnarray*}
and therefore $\phi^{(k)}(m)\in M_n^{(k)}$.

Now apply $\psi^{(k)}$ to $\phi^{(k)}(m)$. Clearly, the highest
exponent of a $y$-variable (say $y_{\pi(n)}$) in $\phi^{(k)}(m)$
is $u = \#\{t\,:\, \la'_t\ge k\}$ $(= \la_k)$. Thus
\begin{eqnarray*}
\psi^{(k)}(\phi^{(k)}(m)) &=&
\prod_{\la'_t\le k-1} (x_{\pi(1)}\cdots x_{\pi(\la'_t)}) \cdot
\prod_{\la'_t\ge k} (x_{\pi(\la'_t+1)}\cdots x_{\pi(n)})^{-1} \cdot
(x_1 \cdots x_n)^u\cr
&=& \prod_{\la'_t\le k-1} (x_{\pi(1)}\cdots x_{\pi(\la'_t)}) \cdot
\prod_{\la'_t\ge k} (x_{\pi(1)}\cdots x_{\pi(\la'_t)}) = m.
\end{eqnarray*}

Conversely, take $m'\in M_n^{(k)}$. We want to show that
$\phi^{(k)}(\psi^{(k)}(m')) = m'$.
We can write
$$
m' = \prod_{i=1}^{k-1} x_{\pi(i)}^{\mu_i} \cdot
     \prod_{i=k+1}^{n} y_{\pi(i)}^{\mu_i},
$$
where $\mu_1\ge \ldots \ge \mu_{k-1}\ge 0$ and $0\le \mu_{k+1}\le
\ldots \le \mu_n$. Here $\pi$ is the unique permutation that
orders first the indices $i\in\supp_x(m')$, then the indices
$i\not\in\supp_x(m') \cup\supp_y(m')$, and then the indices
$i\in\supp_y(m')$. The $x$-indices are ordered by weakly
decreasing exponents, the $y$-indices are ordered by weakly
increasing exponents, and indices with equal exponents are ordered
in increasing (index) order. The variables with a given exponent
are ordered by increasing indices. The highest exponent of a
$y$-variable in $m'$ is $u = \mu_n$, and therefore
$$
\psi^{(k)}(m') = \prod_{i=1}^{k-1} x_{\pi(i)}^{u+\mu_i} \cdot x_{\pi(k)}^u \cdot
                 \prod_{i=k+1}^{n} x_{\pi(i)}^{u-\mu_i} =
                 \prod_{i=1}^{n} x_{\pi(i)}^{\la_i}
$$
where
$$
\la_i = \begin{cases}
u+\mu_i, &\hbox{if $1\le i\le k-1$;}\cr
u, &\hbox{if $i = k$;}\cr
u-\mu_i, &\hbox{if $k+1\le i\le n$.}\cr
\end{cases}
$$
It follows that $\la = (\la_1,\ldots,\la_n)$ is the index partition of
$\psi^{(k)}(m')$, and the conjugate partition $\la'$ satisfies
$$
\la'_t\ge k \iff t\le u = \mu_n.
$$
The map $\phi^{(k)}$ replaces, for each $t$ with $\la'_t\ge k$,
the product $x_{\pi(1)}\cdots x_{\pi(\la'_t)}$
by the product $y_{\pi(\la'_t+1)}\cdots y_{\pi(n)}$.
It therefore reduces by $u$ the exponent of each $x_i$ with $i\le k$,
removes all the $x_i$ with $i\ge k+1$, and replaces them by $y_i$ with
exponents $u-\la_i = \mu_i$. This implies that
$$
\phi^{(k)}(\psi^{(k)}(m')) = \prod_{i=1}^{k-1} x_{\pi(i)}^{\mu_i} \cdot
     \prod_{i=k+1}^{n} y_{\pi(i)}^{\mu_i} = m',
$$
as claimed.

\qed

\medskip

\begin{cor}
The map $\psi^{(k)}: M_n^{(k)}\to M_n^{(n)}$ is a bijection.
\end{cor}

\subsection{Action and Invariants}

\begin{df}
For $1\le m\le n-1$ let
$$
e_m^{(k)}:=\begin{cases} e_m(\bar x) = e_m(x_1,\ldots,x_n), &
\hbox{ if $1\le m\le k-1$}; \cr e_{n-m}(\bar y) =
e_{n-m}(y_1,\ldots,y_n), & \hbox{ if $k\le m\le n-1$},
\end{cases}
$$
where $e_m(\bar x)$ is the $m$-th elementary symmetric function.

\noindent For a partition $\mu=(0 < \mu_1 \leq \cdots \leq \mu_\ell)$
with $\mu_\ell
< n$ let
$$
e_\mu^{(k)}:=\prod\limits_{i=1}^\ell e_{\mu_i}^{(k)}.
$$
For $e_{\mu}^{(n)}$ we shall frequently use the short notation $e_{\mu}$.
\end{df}

Consider the natural $S_n$-action on $\DP_n^{(k)}$, defined by:
$$
\pi(x_i) := x_{\pi(i)},\quad \pi(y_i) := y_{\pi(i)},
\qquad(\forall \pi\in S_n, 1\le i\le n)
$$

\begin{pro}\label{t.action1}
For every $\sigma\in S_n$ and $1\le k\le n$,
$$
\sigma \psi^{(k)}= \psi^{(k)} \sigma .
$$
\end{pro}

\begin{cor}\label{t.action2}
$\DP_n^{(k)}$ and $\DP_n^{(n)}$ are isomorphic as $S_n$-modules.
\end{cor}

Let ${\DP^{(k)}_n}^{S_n}$ be the algebra of $S_n$-invariants in
$\DP^{(k)}_n$. Proposition~\ref{t.action1} implies

\begin{cor}\label{t.invariants1}
For any $1\le k\le n$,
\begin{itemize}
\item[1.] ${\DP^{(k)}_n}^{S_n}$ is generated, as an algebra, by
$e_m^{(k)}$, $1\le m<n$.
\item[2.] The set $\{e_\mu^{(k)}\,:\, \mu = (\mu_1 \leq \cdots
\leq \\mu_\ell) \ \mbox{with} \ \mu_\ell <n\}$ forms a (vector
space) basis for ${\DP^{(k)}_n}^{S_n}$. \item[3.] For every
partition $\mu =(\mu_1 \leq \cdots \mu_\ell)$ with $\mu_\ell < n$,
$$
\psi^{(k)} (e_\mu^{(k)})= e_\mu.
$$
\end{itemize}
\end{cor}

\begin{rem}\label{r.action3}\rm \
By Proposition~\ref{t.action1}, 
$\psi^{(k)}$ and $\phi^{(k)}$ send invariants to invariants.
Unfortunately, these maps are not multiplicative and they do not send
the ideal generated by invariants (with no constant term) to its
analogue. For example, $x_1x_3+x_2x_3+x_3^2=x_3 \cdot
e_1(x_1,x_2,x_3)\in I_{(3)}$ but $\phi^{(2)}(x_1x_3+x_2x_3+x_3^2)=
y_2+y_1+x_3^2\not\in I_{(2,1)}$.
\end{rem}

In order to define the ``correct" map, recall the definitions
of the Garsia-Stanton basis element $a_\pi$ and the $k$-th descent basis element
$a_\pi^{(k)}$ corresponding to $\pi\in S_n$ (see Definition~\ref{d.basis}).
\begin{cla}
For every $\pi\in S_n$ and $1\le k\le n$,
$$
\psi^{(k)}(a_\pi^{(k)})=a_\pi.
$$
\end{cla}

In the next section it will be shown that
$\{a_\pi^{(k)}e_\mu^{(k)}\,:\, \pi\in S_n,\, \mu =(\mu_1 \leq \cdots \mu_\ell) \ \mbox{with} \ \mu_\ell < n\}$ forms a
basis for $\DP^{(k)}_n$. Thus the map
$$
\tilde\psi^{(k)}(a_\pi^{(k)} \cdot e_\mu^{(k)}) :=
\psi^{(k)}(a_\pi^{(k)}) \cdot \psi^{(k)}(e_\mu^{(k)}) =
a_\pi \cdot e_\mu
$$
extends to a linear map $\tilde\psi^{(k)}: \DP^{(k)}_n \to \DP_n^{(n)}$
that clearly sends the $k$-th ideal generated by invariants to its $n$-th analogue.




\section{Straightening}\label{s.straight}

\subsection{Basic Notions}\label{s.st-prel-old}

We now generalize several notions and facts which were used
in~\cite{TAMS} for straightening the coinvariant algebra of type
$A$. Some of them were already introduced in the proof of
Proposition~\ref{t.inverse_maps} above. Proofs are similar to
those in~\cite{TAMS}, and will be omitted.

\medskip

Each monomial $m\in M_n^{(k)}$ can be written in the form
\begin{equation}\label{m*}
m = \prod_{i=1}^{k-1} x_{\pi(i)}^{p_i} \cdot
    \prod_{i=k+1}^{n} y_{\pi(i)}^{p_i},
\end{equation}
where $p_1\ge \ldots \ge p_{k-1}\ge 0$ and $0\le p_{k+1}\le \ldots
\le p_n$. Here $\pi = \pi(m)$, the {\em index permutation} of $m$,
is the unique permutation that orders first the indices
$i\in\supp_x(m)$, then the indices $i\not\in\supp_x(m)
\cup\supp_y(m)$, and then the indices $i\in\supp_y(m)$. The
$x$-indices are ordered by weakly decreasing exponents, the
$y$-indices are ordered by weakly increasing exponents, and
indices with equal exponents are ordered in increasing (index)
order. 

\smallskip

\begin{obs}
The index permutation is preserved by $\psi^{(k)}$, i.e.,
for any monomial $m\in M_n^{(k)}$
$$
\pi(m) = \pi(\psi^{(k)}(m)).
$$
\end{obs}

\begin{cla}\label{t.exp}
Let $m$ be a monomial in $M_n^{(k)}$, $\pi=\pi(m)$ its index permutation,
and $a^{(k)}_\pi$ the corresponding descent basis element:
\begin{eqnarray*}
m &=& \prod_{i=1}^{k-1} x_{\pi(i)}^{p_i}\cdot
      \prod_{i=k+1}^{n} y_{\pi(i)}^{p_i},\cr
a^{(k)}_\pi &=& \prod_{i=1}^{k-1} x_{\pi(i)}^{d^{(k)}_i(\pi)}\cdot
                \prod_{i=k+1}^{n} y_{\pi(i)}^{d^{(k)}_i(\pi)}.
\end{eqnarray*}
Then: the sequence $(p_i-d^{(k)}_i(\pi))_{i=1}^{k-1}$ of
$x$-exponents in $m/a^{(k)}_\pi$ consists of nonnegative integers,
and is weakly decreasing:
$$
p_1-d^{(k)}_1(\pi)\ge \ldots \ge p_{k-1}-d^{(k)}_{k-1}(\pi) \ge 0,
$$
and the sequence $(p_i-d^{(k)}_i(\pi))_{i=k+1}^{n}$ of
$y$-exponents $m/a^{(k)}_\pi$ consists of nonnegative integers,
and is weakly increasing:
$$
0 \le p_{k+1}-d^{(k)}_{k+1}(\pi)\le \ldots \le
p_{n}-d^{(k)}_{n}(\pi).
$$
\end{cla}

\smallskip

For a monomial $m\in M_n^{(k)}$ of the form (\ref{m*}) with index
permutation $\pi\in S_n$,
let the
associated pair of exponent partitions
$$
\la(m) = (\la_x(m),\la_y(m)) :=
((p_1, p_2, \ldots, p_{k-1}),(p_{n},p_{n-1}, \ldots, p_{k+1}))
$$
be its {\em exponent bipartition}. Note that $\la(m)$ is a
bipartition of the total bi-degree of $m$.

Define the {\em complementary bipartition}
$\mu^{(k)}(m)=(\mu_x(m),\mu_y(m))$ of $m$ to be the pair of
partitions conjugate to the partitions
$(p_{i}-d^{(k)}_i(\pi))_{i=1}^{k-1}$ and
$(p_{i}-d^{(k)}_i(\pi))_{i=n}^{k+1}$ respectively; namely,
$$
(\mu_x)_j := |\{1\le i\le k-1\,:\, p_i-d^{(k)}_i(\pi)\ge j\}|
\qquad(\forall j)
$$
and
$$
(\mu_y)_j := |\{k+1\le i\le n\,:\, p_i-d^{(k)}_i(\pi)\ge j\}|
\qquad(\forall j).
$$
If $k=n$ then, for every monomial $m\in M_n^{(n)}$, $\mu_y(m)$ is
the empty partition. In this case we denote
$$
\mu(m) := \mu_x(m).
$$
With each $m\in M_n^{(k)}$ we associate the {\em canonical complementary partition}
$$
\nu(m) := \mu(\psi^{(k)}(m)).
$$

\smallskip

\begin{exa} Let $m = x_1^2 y_2^4 x_3^2 y_5 x_6^3$ with $n=7$ and $k=5$. Then
$$
m = x_6^3 x_1^2 x_3^2 y_5^1 y_2^4,\quad
\la(m) = ((3,2,2,0),(4,1)),\quad
\pi = 6134752\in S_7,
$$
$$
a^{(5)}_\pi=x_6y_5y_2^2,\quad \la(a^{(5)}_{\pi}) =
((1,0,0,0),(2,1)),\quad \mu^{(5)}(m) =
((2,2,2,0)',(2,0)')=((3,3),(1,1)),
$$
$$
\psi^{(5)}(m) = x_6^7 x_1^6 x_3^6 x_4^4 x_7^4 x_5^3,\quad
a_{\pi} = x_6^3 x_1^2 x_3^2 x_4^2 x_7^2 x_5^1,\quad
\nu(m) = \mu(\psi^{(5)}(m)) = (4,4,4,2,2,2)' = (6,6,3,3).
$$
\end{exa}


\smallskip

\begin{df}\label{order}
\begin{itemize}
\item[1.]
For two partitions $\la$ and $\mu$, denote $\la\
\underline{\triangleleft}\ \mu$ if $\la$ is weakly smaller than
$\mu$ in dominance order. For two bipartitions
$\la^1=(\mu^1,\nu^1)$ and $\la^2=(\mu^2,\nu^2)$, denote $\la^1\
\underline{\triangleleft}\ \la^2$ if $\mu^1\
\underline{\triangleleft}\ \mu^2$ and $\nu^1\
\underline{\triangleleft}\ \nu^2$.
\item[2.]
For two monomials
$m_1,m_2\in M^{(k)}_n$ of the same total bi-degree $(p,q)$,
write $m_1\preceq_k m_2$ if:
\begin{itemize}
\item[(1)]
$\la(m_1)\ \underline{\triangleleft} \ \la(m_2)$; and
\item[(2)]
if $\la(m_1) = \la(m_2)$
then
$\inv (\pi(m_1))
> \inv (\pi(m_2))$.
\end{itemize}
\end{itemize}
\end{df}

\subsection{The Straightening Algorithm}

\begin{lem}\label{t.mk}
Let $m\in M^{(k)}_n$ be a monomial. For $1\le d\le n-1$, let
$S_{(d)}$ be the set of all monomials which appear (with
coefficient 1) in the expansion of the polynomial $m\cdot
e_d^{(k)}\in\DP_n^{(k)}$. Let $\pi=\pi(m)$ be the index
permutation of $m$ and denote
$$
m_{(d)} := \begin{cases} m\cdot x_{\pi(1)} \cdots x_{\pi(d)}, &
\hbox{ if $1\le d\le k-1$};\cr m\cdot y_{\pi(d+1)} \cdots
y_{\pi(n)}, & \hbox{ if $k\le d\le n-1$}.
\end{cases}
$$
Then:
\begin{itemize}
\item[(1)] $m_{(d)}\in S_{(d)}$. \item[(2)] $(m'\in S_{(d)} \hbox{
and } m'\ne m_{(d)}) \then m'\prec_k m_{(d)}$.
\end{itemize}
\end{lem}

The proof of Lemma~\ref{t.mk} is similar to the proof
of~\cite[Lemma 3.2]{TAMS} and is omitted.

\begin{cor}\label{t.S}
Let $m\in M^{(k)}_n$ be a monomial, $\pi=\pi(m)$ its index
permutation, and $\nu=\mu(\psi^{(k)}(m))$ the canonical
complementary partition of $m$ defined in
Subsection~\ref{s.st-prel-old}. Let $S$ be the set of monomials
which appear (with nonzero coefficient) in the expansion of
$a_\pi^{(k)}\cdot e_\nu^{(k)}\in\DP_n^{(k)}$. Then :
\begin{itemize}
\item[(1)]
$m\in S$.
\item[(2)]
$(m'\in S \hbox{ and } m'\not=m) \then m'\prec_k m$.
\end{itemize}
\end{cor}

A straightening algorithm follows.

\smallskip

\begin{note}
\noindent{\bf A Straightening Algorithm:}

\noindent For a monomial $m\in \DP_n^{(k)}$, let $\pi=\pi(m)$ be
its index permutation, $a_\pi^{(k)}$ the corresponding descent
basis element, and $\nu=\mu(\psi^{(k)}(m))$ the corresponding
canonical complementary partition. Write (by Corollary \ref{t.S})
$$
m=a_\pi^{(k)}\cdot e_\nu^{(k)} - \Sigma,
$$
where $\Sigma$ is a sum of monomials $m'\prec_k m$. Repeat the
process for each $m'$.
\end{note}

The algorithm implies

\begin{lem}\label{straight.l}{\bf (Straightening Lemma)}
Each monomial $m\in \DP_n^{(k)}$ has an expression
$$
m = a_{\pi(m)}^{(k)}e_{\nu(m)}^{(k)} +
\sum_{m'\prec_k m} n_{m',m} a_{\pi(m')}^{(k)}e_{\nu(m')}^{(k)},
$$
where $n_{m',m}$ are integers.
\end{lem}



The Straightening Lemma yields a direct proof to the following
special case of Theorem~\ref{mk}.

\begin{cor}\label{4.5'}
The set $\{a_\pi^{(k)}\,:\,\pi\in S_n\}$ forms a basis for ${\bf
H}_{(1^{k-1},n-k+1)}'$.
\end{cor}

\noindent{\bf Proof.\ } Recall that $\DP_n^{(k)}:={\bf R}/{\cal
I}_k$, where ${\bf R}:=\bbq[x_1,\dots,x_n,y_1,\dots,y_n]$
(Definition~\ref{3.1}). Let $I^{0}_k$ be the ideal of ${\bf R}$
generated by the elementary symmetric functions $e_d^{(k)}$ $(1\le
d\le n-1)$. By definition, these functions are symmetric functions
of either $\bar x$ or $\bar y$. Thus, Corollary~\ref{s0} and
Claim~\ref{32},
\begin{equation}\label{ideals-inclusion}
{\cal I}_k+I^{0}_k\subseteq J_{(1^{k-1},n-k+1)}.
\end{equation}
By Lemma~\ref{straight.l}, $\{a_\pi^{(k)}: \, \pi\in S_n\}$ spans
${\bf R}/({\cal I}_k+I^{0}_k)$ and therefore also ${\bf
R}/J_{(1^{k-1},n-k+1)}$ as vector spaces over $\bbq$.

On the other hand, by~\cite{GH1},\ $\hbox{dim }({\bf
R}/J_{(1^{k-1},n-k+1)})= |S_n|$, and therefore $\{a_\pi^{(k)}: \,
\pi\in S_n\}$ is actually a basis for ${\bf H}_{(1^{k-1},n-k+1)}'={\bf
R}/J_{(1^{k-1},n-k+1)}$. \qed


\bigskip

\noindent{\bf Second Proof of Theorem~\ref{t.m2}.\ } From the
proof of Corollary~\ref{4.5'} it follows that
$$
\hbox{dim }({\bf R}/({\cal I}_k+I^{0}_k)) = \hbox{dim }({\bf
R}/J_{(1^{k-1},n-k+1)}).
$$
Combining this with inclusion (\ref{ideals-inclusion}), it follows
that
$$
{\cal I}_k+I^{0}_k= J_{(1^{k-1},n-k+1)}.
$$
This completes the proof. \qed



\section{Descent Representations}\label{s.des-rep}

Recall the definition of the Solomon descent representation
$\rho^A$ from Subsection~\ref{m.des-rep}. The following theorem is
well known. Recall the definition of $\Des(T)$ in
Subsection~\ref{Intro-irreducibles}.

\begin{thm}\label{solomon-des-decomp}
For any subset $A\subseteq [n-1]$ and partition $\mu\vdash n$, the
multiplicity in the descent representation $\rho^A$ of the
irreducible $S_n$-representation corresponding to $\mu$ is
$$
m_{\mu}^A := \#\,\{T\in SYT(\mu)\,:\, \Des(T)=A\},
$$
the number of standard Young tableaux of shape $\mu$ with descent
set $A$.
\end{thm}


Recall the definition of a bipartition and the domination order on
bipartitions from Definition~\ref{order}.1.
Recall the subspaces $I_{\la}^{(k)\underline{\triangleleft}}$ and
$I_{\la}^{(k)\triangleleft}$ of the $S_n$-module ${\bf
H}_{(1^{k-1},n-k+1)}'$,   and the quotient $ R_{\la}^{(k)}:=
I_{\la}^{(k)\underline{\triangleleft}}/ I_{\la}^{(k)\triangleleft}
$, defined in Subsection~\ref{m.des-rep}.

\begin{pro}
$I_{\la}^{(k)\underline{\triangleleft}}$,
$I_{\la}^{(k)\triangleleft}$ and thus $R_{\la}^{(k)}$ are
$S_n$-invariant.
\end{pro}

\noindent{\bf Proof.\ } Follows from the straightening algorithm.
\qed

Recall that $\lambda = ((\mu_1 \geq \cdots \geq \mu_k = 0),
(\nu_1 \geq \cdots \geq \nu_{n-k}=0))$ is an $(n,k)$ partition,
if $\forall (1\le i<k-1)\quad(\mu_i - \mu_{i+1}
\in \{0,1\})$ and $\forall  (1\le i<n-k)\quad (\nu_{i} - \nu_{i+1} \in
\{0,1\}).$

\begin{lem}\label{t.basis_R}
Let $\la=(\mu,\nu)$ be an $(n,k)$-bipartition. Then
$R_{\la}^{(k)} \ne \{0\}$ and
$\{a_{\pi}^{(k)}+I_{\la}^{(k)\triangleleft} \,:\, \Des(\pi) =
A_{\la}\}$ is a basis for $R_{\la}^{(k)}$, where $ A_{\la}$ is
defined as in (\ref{e.Rbasis3}).
\end{lem}

The proof is analogous to the proof of~\cite[Corollary
3.10]{TAMS}.

\begin{thm}\label{t.action}
The $S_n$-action on $R_{\la}^{(k)}$ is given by
$$
s_j(a_{\pi}^{(k)}) = \begin{cases} a_{s_j\pi}^{(k)}, &\mathrm{if\
} |\pi^{-1}(j+1) - \pi^{-1}(j)| > 1;\cr a_{\pi}^{(k)},
&\mathrm{if\ } \pi^{-1}(j+1) = \pi^{-1}(j)+1;\cr
-a_{\pi}^{(k)}-\sum_{\sigma\in A_j(\pi)} a_{\sigma}^{(k)},
  &\mathrm{if\ } \pi^{-1}(j+1) = \pi^{-1}(j)-1.\cr
\end{cases}
$$
Here $s_j = (j,j+1)$ $(1\le j< n$) are the Coxeter generators of
$S_n$; $\{a_{\pi}^{(k)}+I_{\la}^{(k)\triangleleft} \,:\,
\Des(\pi)=A_{\la}\}$
is the descent basis of $R_{\la}^{(k)}$; for $\pi\in
S_{\la}$ with $\pi^{-1}(j+1) = \pi^{-1}(j)-1$ we define
\begin{eqnarray*}
t &:=& \pi^{-1}(j+1),\cr m_1 &:=& \max\{i\in\Des(\pi)\cup\{0\}
\,:\, i\le t-1\},\cr m_2 &:=& \min\{i\in\Des(\pi)\cup\{n\} \,:\,
i\ge t+1\}
\end{eqnarray*}
(so that $\pi(t) = j+1$, $\pi(t+1) = j$, and
$\{m_1+1,\ldots,m_2\}$ is the maximal interval containing $t$ and
$t+1$ on which $s_j\pi$ is increasing); and let $A_j(\pi)$ be the
set of all $\sigma\in S_n$ satisfying
\begin{enumerate}
\item $(i\le m_1 \mathrm{\ or\ } i\ge m_2+1) \,\then\, \sigma(i) =
\pi(i);$ \item the sequences $(\sigma(m_1+1),\ldots,\sigma(t))$
and $(\sigma(t+1),\ldots,\sigma(m_2))$ are increasing; \item
$\sigma \not\in\{\pi, s_j\pi\}$ (i.e., $\{\sigma(t),\sigma(t+1)\}
\ne \{j,j+1\}$).
\end{enumerate}
\end{thm}

\begin{exa}
Let $\pi = 2416573\in S_7$ and $j=5$. Then:
$$
j=5,\,j+1=6;\,\, t=4,\,t+1=5;
$$
$$
\Des(\pi) =\{2,4,6\};\,\, m_1=2,\,m_2=6;\,\, s_j\pi =
24\underline{1567}3;
$$
$$
A_j(\pi) = \{24\underline{1756}3,\, 24\underline{5617}3,\,
             24\underline{5716}3,\, 24\underline{6715}3\}.
$$
Note that $|A_j(\pi)| = {m_2 - m_1 \choose t - m_1} - 2 = {4
\choose 2} - 2 = 4$.
\end{exa}

\begin{cor}
The $S_n$-representation on $R_{\la}^{(k)}$ is independent of $k$.
\end{cor}

\noindent {\bf Proof of Theorem~\ref{t.action}.\ } If
$|\pi^{-1}(j+1) - \pi^{-1}(j)| > 1$ (i.e., if $j$ and $j+1$ are
not adjacent in the sequence $(\pi(1),\ldots,\pi(n))$), then
$\Des(s_j\pi) = \Des(\pi)$ and therefore $s_j(a_{\pi}^{(k)}) =
a_{s_j\pi}^{(k)}$.

If $\pi^{-1}(j+1) = \pi^{-1}(j)+1$ then $j$ immediately precedes
$j+1$ in the sequence $(\pi(1),\ldots,\pi(n))$. Therefore $t :=
\pi^{-1}(j) \not\in\Des(\pi)$, so that: if $t+1\le k-1$ then $x_j$
and $x_{j+1}$ have the same exponent in $a_{\pi}^{(k)}$; if $t\ge
k+1$ then $y_j$ and $y_{j+1}$ have the same exponent in
$a_{\pi}^{(k)}$; and if $t\in\{k-1,k\}$ then none of the variables
$x_j$, $x_{j+1}$, $y_j$ and $y_{j+1}$ appears in $a_{\pi}^{(k)}$.
In all of these cases we have $s_j(a_{\pi}^{(k)}) =
a_{\pi}^{(k)}$.

We are left with the most involved case: $\pi^{-1}(j+1) =
\pi^{-1}(j)-1$. Recall the notations $t$, $m_1$, $m_2$ and
$A_j(\pi)$ from the statement of the theorem. In particular,
recall that $\pi(t) = j+1$, $\pi(t+1) = j$, and
$\{m_1+1,\ldots,m_2\}$ is the maximal interval containing $t$ and
$t+1$ on which $s_j\pi$ is increasing.

Assume first that $m_2\le k$. The variables
$x_{\pi(m_1+1)},\ldots,x_{\pi(m_2)}$ are all the $x$-variables
having the same exponent in $a_{s_j\pi}^{(k)}$ as $x_j$ (and
$x_{j+1}$). Since $t\in\Des(\pi)$, the exponents of
$x_{\pi(1)},\ldots,x_{\pi(t)}$ in $a_{\pi}^{(k)}$ are $1$ higher
than the corresponding exponents in $a_{s_j\pi}^{(k)}$. Thus the
product $a_{s_j\pi}^{(k)} \cdot e_t(x)$ (in $\DP_n^{(k)}$) is a
sum of $n \choose t$ monomials, one of which is $a_{\pi}^{(k)}$.
In all of these monomials, $t$ of the $x$-variables have their
exponent increased by $1$ (with respect to $a_{s_j\pi}^{(k)}$). If
these $t$ variables miss any one of
$x_{\pi(1)},\ldots,x_{\pi(m_1)}$, or include any one of
$x_{\pi(m_2+1)},\ldots,x_{\pi(n)}$, then the monomial belongs to
$J_{\la}^{(k)\triangleleft}$ and contributes nothing to
$R_{\la}^{(k)}$. The remaining ${m_2 - m_1 \choose t - m_1}$
monomials are: $a_{\sigma}^{(k)}$ for $\sigma\in A_j(\pi)$,
$a_{\pi}^{(k)} = a_{s_j\pi}^{(k)} \cdot x_{\pi(1)}\cdots
x_{\pi(t)}$, and $s_j(a_{\pi}^{(k)}) = a_{s_j\pi}^{(k)} \cdot
x_{s_j\pi(1)}\cdots x_{s_j\pi(t)}$. On the other hand,
$a_{s_j\pi}^{(k)} \cdot e_t(x) \equiv 0$ in $R_{\la}^{(k)}$ since
$1\le t\le m_2-1\le k-1$. This proves the claim of the theorem in
this case.

If $m_1+1\ge k$, then an analogous argument holds for $y$- instead
of $x$-variables, and here $a_{s_j\pi}^{(k)} \cdot e_{n-t}(y)
\equiv 0$ in $R_{\la}^{(k)}$ since $1\le n-t\le n-m_1-1\le n-k$.

Finally, if $m_1+1< k< m_2$, then the variables
$x_{\pi(m_1+1)},\ldots,x_{\pi(m_2)}$ (and the corresponding
$y$-variables) do not appear at all in $a_{s_j\pi}^{(k)}$. If
$t\le k-1$ (respectively, $t\ge k$) then $a_{s_j\pi}^{(k)} \cdot
e_t(x)$ (respectively, $a_{s_j\pi}^{(k)} \cdot e_{n-t}(y)$) is,
again, a sum of ${m_2 - m_1 \choose t - m_1}$ monomials:
$a_{\sigma}^{(k)}$ for $\sigma\in A_j(\pi)$, $a_{\pi}^{(k)} =
a_{s_j\pi}^{(k)} \cdot x_{\pi(1)}\cdots x_{\pi(t)}$ (respectively,
$a_{\pi}^{(k)} = a_{s_j\pi}^{(k)} \cdot y_{\pi(t+1)}\cdots
y_{\pi(n)}$), and $s_j(a_{\pi}^{(k)}) = a_{s_j\pi}^{(k)} \cdot
x_{s_j\pi(1)}\cdots x_{s_j\pi(t)}$ (respectively,
$s_j(a_{\pi}^{(k)}) = a_{s_j\pi}^{(k)} \cdot y_{s_j\pi(t+1)}\cdots
y_{s_j\pi(n)}$). This completes the proof.

\qed

\bigskip


\begin{thm}\label{t.SO}
Let $\la=(\mu,\nu)$ be an $(n,k)$-bipartition. $R_{\la}^{(k)}$ is isomorphic,
as an $S_n$-module, to the Solomon descent representation determined by
the descent class $\{\pi\in S_n\,:\,\Des(\pi)=A_\la\}$.
\end{thm}

\noindent{\bf Proof.} By Theorem~\ref{t.action} together with
Lemma~\ref{t.basis_R}, for every Coxeter generator $s_i$, the
representation matrices of $s_i$ on $R_{\la}^{(k)}$  and on
$R_\la^{(n)}$ with respect to the corresponding $k$-th and $n$-th
descent monomials respectively are identical. By~\cite[Theorem 4.1]{TAMS},
the multiplicity of the irreducible $S_n$-representation corresponding
to $\mu$ in $R_\la^{(n)}$ is
$m_{S,\mu} := \#\,\{\,T\in SYT(\mu)\, :\, \Des(T)=A_\la\,\}$, the
number of standard Young tableaux of shape $\mu$ and descent set
$A_\la$. Theorem~\ref{solomon-des-decomp} completes the proof.

\qed

\bigskip

Let $R^{(k)}_{t_1,t_2}$ be the $(t_1,t_2)$-th homogeneous
component of ${\bf H}_{(1^{k-1},n-k+1)}'$.

\begin{cor}\label{zigzags}
For every   $t_1, t_2\ge 0$ and $1\le k\le n$, the
$(t_1,t_2)$-th homogeneous component of ${\bf H}_{(1^{k-1},n-k+1)}'$
decomposes into a direct sum of Solomon descent representations
as follows:
$$
R^{(k)}_{t_1,t_2} \cong \bigoplus_{\la} R_{\la}^{(k)},
$$
where the sum is over all $(n,k)$-bipartitions and
$$
\sum\limits_{\nu_i>\nu_{i+1}{\rm\ and\ } i\ge k} (n-i) =t_1, \qquad
\sum\limits_{\mu_i>\mu_{i+1}{\rm\ and\ } i<k} i =t_2
$$
\end{cor}

Finally, Theorem~\ref{t.SO} implies Stembridge's Theorem~\ref{t.stemb}.

\medskip

\noindent{\bf First Proof of Theorem~\ref{t.stemb}.} Combine
Theorems~\ref{solomon-des-decomp} and~\ref{t.SO} with
Corollary~\ref{zigzags}.

\qed

\section{The Schur Function Expansion of
$\tilde{H}_{(1^{k-1},n-k+1)}(\ox;q,t)$}\label{s.decomposition}

In this section we give a direct combinatorial proof of Theorems~\ref{t.stemb}
and~\ref{t.stemb2}, using the axiomatic characterization of Macdonald polynomials
and properties of the RSK algorithm.

\subsection{Preliminaries}\label{s.prelim}

A {\em skew} Young diagram $\la/\mu$ is the difference of Young diagrams
$\la$ and $\mu \subseteq \la$. A skew diagram is a {\em horizontal strip}
(resp. {\em vertical strip}) if it does not contain two cells in the same
 column (row). A {\em semistandard Young tableau} of (skew) shape $\la$ is
a function $T$ from the diagram of $\la$ to the ordered alphabet
$${\it A}_+ = \{ 1 < 2 < \ldots \}$$
which is weakly increasing in each row of $\la$ from left to right
and strictly increasing in each column from bottom to top. A
semistandard tableau is {\em standard} if it is a bijection form
$\la$ to $\{1, \ldots, n = |\la|\}$. More generally, we shall
admit the ordered alphabet
$${\it A}_{\pm} = {\it A}_+ \cup {\it A}_- = \{1 < \bar{1} < 2 < \bar{2} < \cdots \}$$
of {\em positive} letters $1,2, \ldots $ and {\em negative}
letters $\bar{1},\bar{2}, \ldots $. For any alphabet $B$, we let
$B^*$ denote the set of all words over $B$ and $B^n$ denote the
set of all words $w \in B^*$ of length $n$.

A {\em super} tableau is a function $T:\la \rightarrow {\it
A}_\pm$ which is weakly increasing in each row and column and is
such that the entries equal to $a$ occupy a horizontal strip if
$a$ is positive and a vertical strip if $a$ is negative. Thus a
semistandard tableaux is just a super tableau where all the
entries are positive. We denote
\begin{eqnarray*}
SSYT(\la) &=& \{\mbox{semistandard tableaux} \ T:\la \rightarrow {\it A}_+\}\\
SSYT_\pm (\la) &=& \{\mbox{super tableaux} \ T:\la \rightarrow {\it A}_\pm\}\\
SSYT(\la,\mu) &=& \{\mbox{semistandard tableaux} \ T:\la\rightarrow {\it A}_+ \
\mbox{with entries} \ 1^{\mu_1}, 2^{\mu_2}, \ldots \}\\
SSYT_\pm (\la,\mu,\nu) &=&
\{\mbox{semistandard  super tableaux} \ T:\la \rightarrow {\it A}_+ \
\mbox{with entries} \ 1^{\mu_1}, \bar{1}^{\nu_1}, 2^{\mu_2}, \bar{2}^{\nu_2}, \ldots\} \\
SYT(\la) &=& SSYT(\la,1^n)=
\{\mbox{standard tableaux} \ T:\la \rightarrow \{1,\ldots,n = |\la|\}\}
\end{eqnarray*}
If $T$ is any one of these types of tableaux, we shall let $sh(T) = \la$
denote the shape of $T$, $pos(T)$ denote the number of
positive letters in the range of $T$ and $neg(T)$ denote the
number of negative letters in the range of $T$.

We write $\langle\ ,\ \rangle$ for the Hall inner product on symmetric functions,
defined by either one of the identities
\begin{equation}\label{Hscalar}
\langle h_\la,m_\mu \rangle = \delta_{\la,\mu} = \langle s_\la,
s_\mu \rangle.
\end{equation}
We denote by $\omega$ the involution defined by either one of the identities
\begin{equation}\label{omega}
\omega(h_\la) = e_\la, \ \omega(e_\la) = h_\la, \ \omega(s_\la) =
s_{\la'}
\end{equation}

We shall use square brackets $f[A]$ to denote the plethystic
evaluation of a symmetric function $f$ at a polynomial, rational
function, or formal power series $A$. This is defined by writing $f$ in
terms of the power symmetric functions and then substituting
$p_m[A]$ for $p_m(\ox)$, where $p_m[A]$ is the result of
substituting $a \rightarrow a^m$ for every indeterminate $a$. The
standard $\lambda$-ring identities hold for plethystic evaluation,
e.g. $s_\la[X+Y] = \sum_{\mu \subseteq \la} s_\mu[X]
s_{\la/\mu}[Y]$, etc. In particular, setting $Z = z_1 + z_2 +
\cdots $, we have $f[Z] = f(z_1, z_2, \ldots )$.

If $W= w_1 + w_2 + \ldots$  and $f$ is a symmetric function, we
shall use the notation
\begin{equation}\label{super1}
\omega^W f[Z+W]
\end{equation}
to denote the result of applying $\omega$ to $f[Z+W]= f(z_1,z_2,
\ldots, w_1,w_2 \ldots )$ considered as a symmetric function in
the $w$ variables with functions of $z$ as coefficients. Equations
(\ref{omega}) and (\ref{super1})  then imply that the coefficient
of a monomial $z^{\mu}w^{\eta} =z_1^{\mu_1} z_2^{\mu_2} \cdots
w_1^{\eta_1} w_2^{\eta_2} \cdots $ in $\omega^W f[Z+W]$ is given
by
\begin{equation}
\omega^W f[Z+W]|_{z^{\mu} w^{\eta}} = \langle f, e_{\eta}(z)
h_\mu(w) \rangle.
\end{equation}

If $T$ is a semistandard tableau of (skew) shape $\la$, we set
\begin{equation}
z^T = \prod_{x \in \la} z_{T(x)}.
\end{equation}
Then the usual combinatorial definition of the Schur function
$s_\la$ is
\begin{equation}\label{schur}
s_\la(z_1, z_2, \ldots ) = \sum_{T \in SSYT(\la)}z^T.
\end{equation}
Throughout what follows, we fix
$$Z = z_1 + z_2 + \cdots, \ W = w_1 + w_2 + \cdots,$$
and make the convention that
$$z_{\bar{a}} \ stands \ for \ w_a, \ for \ every \ negative \ letter \ \bar{a} \in {\it A}_- .$$
Then the ``super'' analogue of  (\ref{schur}) is
\begin{equation}\label{hschur}
HS_\la(z,w) = \omega^W s_\la[Z+W] = \sum_{T \in SSYT_\pm (\la)}z^T
\end{equation}
We shall refer to $HS_\lambda(z,w)$ as the super Schur function,
but it is also called a hook Schur function $HS_\la(z,w)$ by
Berele and Regev \cite{BReg87}. We note that the  right hand side
of (\ref{hschur}) is actually independent of the relative order of
positive letters ${\it A}_+$ and negative letters ${\it A}_{-}$.
That is, if ${\it A}_+$ is ordered $1 < 2 < \cdots$ and
${\it A}_-$ is ordered $\bar{1} < \bar{2} < $, then the RHS of
(\ref{hschur}) is independent of the relative order between the
positive and negative letters. Thus, for example, if we let $1 < 2
< \cdots < \bar{1} < \bar{2} < \cdots$, the a super tableaux $T$
of shape $\la$ will consists of a pair of tableau $(T_1,T_2)$
where $T_1:\mu \rightarrow {\it A}_+$ is a semistandard tableau of
shape $\mu$ for some $\mu \subseteq \la$ and $T_2:\la'/\mu'
\rightarrow {\it A}_-$ is a semi-standard tableau of shape
$\la'/\mu'$ or, equivalently, we can think of $T_2$ as a row
strict tableau of shape $\la/\mu$ by conjugation. Thus we can
think of $T$ as a filling of $\la$ with positive and negative
numbers such that positive numbers form a column strict tableau of
shape $\mu$ and the negative numbers form a row strict tableau of
shape $\la/\mu$. Given such a $T$, we can achieve a super tableau
$T'$ from $T$ corresponding to any other relative order between
${\it A}_+$ and ${\it A}_-$ by using jeu de taquin to move the
negative numbers past the positive numbers as was done in
\cite{Rem84,Rem90}. Another proof of the independence of the RHS
of (\ref{schur}) of the relative order of the positive and
negative letters can be found in \cite{RS02}. It follows that
\begin{equation}\label{hschur2}
HS_\la(z,w) = \sum_{\mu \subseteq \la} s_\mu [Z] s_{\la'/\mu'}
[W].
\end{equation}
Similarly, it follows that the RHS of (\ref{hschur}) is also
independent of the relative order of the positive letters among
themselves and the relative order of the negative letters among
themselves.

\subsection{Second Proof of Theorem~\ref{t.stemb}}

As noted in \cite{HHL05b}, the set $\{\tilde{H}_\mu[Z;q,t]\}_{\mu\vdash n}$
can be characterized as the unique basis of the space of homogeneous symmetric polynomials
of degree $n$ over the field $F=\mathbf{Q}(q,t)$ of rational functions
in $q$ and $t$, that satisfies the following three properties:

\begin{description}
\item[A1]
$\tilde{H}_\mu[Z(q-1);q,t] = \sum_{\rho \underline{\triangleleft} \mu'}
c_{\rho,\mu} m_\rho[Z]$.

\item[A2]
$\tilde{H}_\mu[Z(t-1);q,t] = \sum_{\rho \underline{\triangleleft} \mu}
d_{\rho,\mu} m_\rho[Z]$.

\item[A3]
$\tilde{H}_\mu[z_1,\ldots,z_n;q,t]|_{z_1^n} = 1$.

\end{description}
where $\underline{\triangleleft}$ denotes the dominance order.

\bigskip

The row insertion algorithm of the Robinson-Schensted-Knuth (RSK)
correspondence has the property that if $a_1 \cdots a_n$ is a
sequence of positive letters and $$a_1\cdots a_n \rightarrow_{RSK}
(P,Q)$$ where $P$ is a semistandard Young tableau, $Q$ is a standard
tableau, $sh(P)= sh(Q)$, and
\begin{enumerate}
\item if $a_i \leq a_{i+1}$, then $i+1$ lies strictly to right and
weakly below $i$ in $Q$ and \item if $a_i > a_{i+1}$, then $i+1$
lies strictly above and weakly to the left of $i$ in $Q$.
\end{enumerate}
Given a sequence of positive letters $a = a_1 \ldots a_n$, some $1\le k\le n$,
and some linear order $\prec$ on ${\it A}_+$, we let
(in analogy with the corresponding definitions for tableaux
in Subsection~\ref{Intro-irreducibles} above)
\begin{equation}\label{wordmajk}
\maj_{1,k}(a,\prec) = \sum_{1 \leq i < k, a_i \succ a_{i+1}} i
\end{equation}
and
\begin{equation}\label{wordcomajk}
\rmaj_{k,n}(a,\prec) = \sum_{k \leq i < n, a_i \succ a_{i+1}} (n-i)
\end{equation}
It then follows from Theorem \ref{t.stemb} and the properties of
the RSK algorithm that
\begin{eqnarray}\label{Word1}
\tilde{H}_{(1^{k-1},n-k+1)}[Z;q,t]
&=& \sum_{\la \vdash n} s_\la[Z] \tilde{K}_{\la,(1^{k-1},n-k+1)}(q,t) \\
&=& \sum_{\la \vdash n} \sum_{P \in SSYT(\la)} z^{P} \sum_{Q \in SYT(\la)}
q^{\maj_{1,n-k+1}(Q)}t^{\rmaj_{n-k+1,n}(Q)} \nonumber \\
&=& \sum_{a= a_1 \cdots a_n \in {\it A}_{\pm}^n} z_{a_1}\cdots z_{a_n}
q^{\maj_{1,n-k+1}(a,\prec)} t^{\rmaj_{n-k+1,n}(a,\prec)}.
\nonumber
\end{eqnarray}

To prove (A1) and (A2), we need to interpret
$\tilde{H}_{\mu}[Z(q-1)]$ and $\tilde{H}_{\mu}[Z(t-1)]$. Note that
\begin{eqnarray}\label{Word2}
\tilde{H}_{\mu}[Z(q-1)] &=& \sum_{\la \vdash |\mu|}
s_\la[qZ-Z] \tilde{K}_{\la,\mu}(q,t) \\
&=& \sum_{\la \vdash |\mu|} \tilde{K}_{\la,\mu}(q,t) \sum_{\nu
\subseteq \la} s_\nu[qZ] (-1)^{|\la/\nu|} s_{\la'/\nu'}[Z]
\nonumber
\end{eqnarray}
Here we have used the $\la$-ring identity
$$
s_\la[X-Y] = \sum_{\nu \subseteq \la} s_\nu[X] (-1)^{|\la/\nu|}
s_{\la'/\nu'}[Y].
$$
Our goal is to get an interpretation of
$\tilde{H}_{(1^k,n-k)}[Z(1-q)]$ in terms of statistics on words
which is similar to (\ref{Word1}). To this end, we shall consider
an extension of the RSK algorithm to words over ${\it A}_\pm$
where we row insert positive letters and dual row insert negative
letters as in \cite{Rem84}. Recall that in the dual row insertion
algorithm of the Robinson-Schensted-Knuth (DRSK) correspondence,
an $x$ bumps the first element which is greater than or equal to
$x$ in a row into which $x$ is inserted as opposed to $x$ bumping
the first element which is greater than $x$ in a row into which
$x$ inserted in row insertion algorithm. The DRSK algorithm has
the property that that if $\bar{a}_1 \cdots \bar{a}_n$ is a
sequence of negative letters and
$$\bar{a}_1\cdots \bar{a}_n \rightarrow_{DRSK} (P,Q)$$
where the transpose of $P$, $P^T$, is a semistandard Young tableau, $Q$ is a standard
tableau, and $sh(P)= sh(Q)$, then
\begin{enumerate}
\item if $\bar{a}_i <  \bar{a}_{i+1}$, then $i+1$ lies strictly to
right and weakly below $i$ in $Q$ and \item if $\bar{a}_i \geq
\bar{a}_{i+1}$, then $i+1$ lies strictly above and weakly to the
left of $i$ in $Q$.
\end{enumerate}

Given a word $w= w_1 \cdots w_n \in {\it A}_\pm^n$, $1\le k \le n$,
and some fixed linear order $\prec$ on ${\it A}_{\pm}$, we shall let
\begin{enumerate}
\item
$neg(w)$ be the number of negative letters in $w$,
\item
$pos(w)$ be the number of positive letters in $w$,
\item
$Des(w,\prec)$ be the set of all $i$ such that
\begin{description}
\item[(i)]
$a_i \succeq a_{i+1}$ if both $a_i$ and $a_{i+1}$ are negative and
\item[(ii)]
$a_i \succ a_{i+1}$ otherwise;
\end{description}
\item
$\maj_{1,k}(w,\prec) = \sum_{1 \leq i < k, i\in Des(a,\prec)} i$, and
\item
$\rmaj_{k,n} = \sum_{k \leq i < n, i \in Des(a,\prec)} (n-i)$.
\end{enumerate}

This given, to prove Stembridge's Theorem~\ref{t.stemb} directly, we shall
consider the following family of polynomials
$\{\overline{H}_{\mu}[Z;q,t]\}_{\mu \vdash n}$ where
\begin{description}
\item[(a)]
$\overline{H}_{\mu}[Z;q,t] = \tilde{H}_{\mu}[Z;q,t]$ if $\mu \neq (1^{k-1},
n-k+1)$ and
\item[(b)]
\begin{eqnarray}\label{newdef}
\overline{H}_{(1^{k-1},n-k+1)}[Z;q,t]
&=& \sum_{a= a_1 \cdots a_n \in {\it A}_{+}^n} z_{a_1}\cdots z_{a_n}
q^{\maj_{1,n-k+1}(a,\prec)}t^{\rmaj_{n-k+1,n}(a,\prec)} \\
&=& \sum_{\la \vdash n} \sum_{P \in SSYT(\la)}\sum_{Q \in SYT(\la)}
z^{P} q^{\maj_{1,n-k+1}(Q)} t^{\rmaj_{n-k+1,n}(Q)} \nonumber \\
&=& \sum_{\la \vdash n} s_\la [Z] \sum_{Q \in SYT(\la)}
q^{\maj_{1,n-k+1}(Q)} t^{\rmaj_{n-k+1,n}(Q)}.\nonumber
\end{eqnarray}
\end{description}
To prove that $\overline{H}_{\mu}[Z;q,t] = \tilde{H}_{\mu}[Z;q,t]$
for all $\mu$, we need only prove that the family of polynomials
$\{\overline{H}_{\mu}[Z;q,t]\}_{\mu \vdash n}$ satisfies the
analogues of (A1), (A2), and (A3). Clearly from our definitions,
we need only show that $\overline{H}_{(1^{k-1},n-k+1)}[Z;q,t]$ satisfies
the analogues of (A1), (A2), and (A3).

It immediately follows from (\ref{newdef}) that the analogue of
(A3) holds since the only word contributing to the coefficient of
$z_1^n$ is the word $a =1^n$ and clearly
${\maj_{1,n-k+1}(a,\prec)}=
{\rmaj_{n-k+1,n}(a,\prec)} =0$ for all $\prec$ in this
case. Thus $\overline{H}_{(1^{k-1},n-k+1)}[Z;q,t]|_{z_1^n} = 1$ as
desired.

It follows from the properties of the row insertion and dual row
insertion of the RSK correspondence described above that
\begin{eqnarray}\label{supnewdefq}
\overline{H}_{(1^k,n-k)}[Z(q-1);q,t]
&=& \sum_{\la \vdash n}
\left(\sum_{Q \in SYT(\la)} q^{\maj_{1,n-k+1}(Q)}t^{\rmaj_{n-k+1,n}(Q)}\right)
\sum_{\nu \subseteq \la} q^{|\nu|} s_\nu[Z] (-1)^{|\la/\nu|} s_{\la'/\nu'}[Z]
\nonumber \\
&=& \sum_{\la \vdash n} \sum_{P \in SSYT_{\pm}(\la)} \sum_{Q \in SYT(\la)}
(-1)^{neg(P)} q^{pos(P)} z^P q^{\maj_{1,n-k+1}(Q)}t^{\rmaj_{n-k+1,n}(Q)}
\nonumber \\
&=& \sum_{a = a_1 \cdots a_n \in {\it A}_{\pm}^n}
q^{pos(a)}(-1)^{neg(a)} z_{|a_1|}\cdots z_{|a_n|}
q^{\maj_{1,n-k+1}(a,\prec)}t^{\rmaj_{n-k+1,n}(a,\prec)}
\nonumber
\end{eqnarray}
for any linear order $\prec$ on ${\it A}_\pm$. Similarly,
\begin{eqnarray}
\overline{H}_{(1^k,n-k)}[Z(t-1);q,t]
&=& \sum_{\la \vdash n}
\left( \sum_{Q \in SYT(\la)} q^{\maj_{1,n-k+1}(Q)}t^{\rmaj_{n-k+1,n}(Q)}\right)
s_\la[tZ-Z] \\
&=& \sum_{\la \vdash n} \left( \sum_{Q \in SYT(\la)}
q^{\maj_{1,n-k+1}(Q)}t^{\rmaj_{n-k+1,n}(Q)}\right)
\sum_{\nu \subseteq \la} t^{|\nu|} s_\nu[Z] (-1)^{|\la/\nu|} s_{\la'/\nu'}[Z]
\nonumber \\
&=& \sum_{\la \vdash n} \sum_{P \in SSYT_{\pm}(\la)} \sum_{Q \in SYT(\la)}
(-1)^{neg(P)} t^{pos(P)} z^P q^{\maj_{1,n-k+1}(Q)}t^{\rmaj_{n-k+1,n}(Q)}
\nonumber \\
&=& \sum_{a = a_1 \cdots a_n \in {\it A}_{\pm}^n}
t^{pos(a)}(-1)^{neg(a)} z_{|a_1|}\cdots z_{|a_n|}
q^{\maj_{1,n-k+1}(a,\prec)}t^{\rmaj_{n-k+1,n}(a,\prec)}
\nonumber
\end{eqnarray}
for any linear order $\prec$ on ${\it A}_\pm$.

Thus to prove (A1), we must prove that there exists a linear order
$\prec$ on ${\it A}_\pm$ such that
\begin{eqnarray}\label{AA1}
&&\sum_{a = a_1 \cdots a_n \in {\it A}_{\pm}^n}
q^{pos(a)}(-1)^{neg(a)} z_{|a_1|}\cdots z_{|a_n|}
q^{\maj_{1,n-k+1}(a,\prec)}t^{\rmaj_{n-k+1,n}(a,\prec)}= \nonumber \\
&& \sum_{\rho \underline{\triangleleft} (1^{n-k},k)} c_{\rho,(1^{k-1},{n-k+1})}m_\rho.
\end{eqnarray}

Similarly to prove (A2), we must prove that there exists a linear
order $\prec$ on ${\it A}_\pm$ such that
\begin{eqnarray}\label{AA2}
&&\sum_{a = a_1 \cdots a_n \in {\it A}_{\pm}^n}
t^{pos(a)}(-1)^{neg(a)} z_{|a_1|}\cdots z_{|a_n|}
q^{\maj_{1,n-k+1}(a,\prec)}t^{\rmaj_{n-k+1,n}(a,\prec)}
= \nonumber \\
&&\sum_{\rho \underline{\triangleleft} (1^{k-1},n-k+1)} d_{\rho,(1^{k-1},{n-k+1})}m_\rho
\end{eqnarray}
for some $c_{\rho,(1^{k-1},{n-k+1})}, d_{\rho,(1^{k-1},{n-k+1})} \in Q(q,t)$.

We shall prove (\ref{AA1}) and (\ref{AA2}) via simple involutions.
We shall start with proving (\ref{AA1}). For any ordering $\prec$,
define the weight $U(a)$ of a word $a = a_1 \cdots a_n \in {\it
A}_{\pm}^n$ by setting
$$U(a) = q^{pos(a)}(-1)^{neg(a)} z_{|a_1|}\cdots z_{|a_n|}
q^{\maj_{1,n-k+1}(a,\prec)}t^{\rmaj_{n-k+1,n}(a,\prec)}.$$
First we define $\prec$ so that
$$
1 \prec 2 \prec \cdots \prec n \prec \overline{n} \prec
\overline{n-1} \prec \cdots \prec \overline{1}.
$$
Then we define an involution $I_k$ on ${\it A}_{\pm}^n$ as follows
. We let $|i| = |\overline{i}| =i$. Given a word $w = a_1 \cdots
a_n \in {\it A}_{\pm}^n$, we look for the smallest letter $j$
which is repeated in $|a_1| \cdots |a_{n-k+1}|$. If there is no such
letter, we let $I_k(w) =w$. Otherwise, let $a_t$ be right most
occurrence of either $j$ or $\overline{j}$ in $a_1 \cdots
a_{n-k+1}$. Thus $1 < t \leq n-k+1$. Now let $i$ be the smallest
letter in $|a_1| \cdots |a_{t-1}|$ and let $a_s$ be the left most
occurrence of $i$ or $\overline{i}$ in $a_1 \ldots a_{t-1}$. Note
that $1 \leq s \leq n-k$ since $t \leq n-k+1$. Then $I_k(a) = b= b_1
\cdots b_n$ where (i) $b_r = a_r$ if $r \neq s$, (ii) $b_s = i$ if
$a_s = \overline{i}$,and (iii) $b_s = \overline{i}$ if $a_s = i$.
In other words, $I_k(a)$ is the result of changing $a_s$ to
$\overline{i}$ if $a_s =i$ or changing $a_s$ to $i$ if $a_s =
\overline{i}$. Clearly $I_k^2$ is the identity. We claim that
$U(a) = - U(b)$. Clearly $(-1)^{neg(a)} = -(-1)^{neg(b)}$ since we
changed the sign of one letter. We did not change the absolute
value of any letter so that $z_{|a_1|}\cdots z_{|a_n|} =
z_{|b_1|}\cdots z_{|b_n|}$. Since we did not change any of the
letters after place $n-k$, we have
$t^{\rmaj_{n-k+1,n}(a,\prec)} = t^{\rmaj_{n-k+1,n}(b,\prec)}$.
Thus we need only show that
\begin{equation}
q^{pos(a)}q^{\maj_{1,n-k+1}(a,\prec)} =
q^{pos(b)}q^{\maj_{1,n-k+1}(b,\prec)}
\end{equation}
Without loss of generality, we may assume that $a_s =i$. Then there are two cases.\\
{\bf Case 1.} $s =1$. Since $i$ is the smallest element in $|a_1|
\cdots |a_{t}|$, it follows $i \preceq a_2$ and hence $1 \notin
Des(a,\prec)$. We claim that $1 \in Des(b,\prec)$.  Now if $|a_2|
>i$ or $a_2 =i$, then clearly $\overline{i} \succ a_2$ by our
choice of the order $\prec$ so that $1 \in Des(b,\prec)$.  Finally
if $a_2 = \overline{i}$, then $1 \in Des(b,\prec)$ by our
definition of $Des(a,\prec)$ since we have two consecutive equal
negative numbers. Thus it follows that
$\maj_{1,n-k+1}(b,\prec)-1= \maj_{1,n-k+1}(a,\prec)$
and hence
$$q^{pos(a)}q^{\maj_{1,n-k+1}(a,\prec)} = q^{pos(b)+1}q^{\maj_{1,n-k+1}(b,\prec)-1}
=q^{pos(b)}q^{\maj_{1,n-k+1}(b,\prec)}
$$
as desired.\\
\ \\
{\bf Case 2.} $s >1$.  We can argue exactly as in Case 1 that $s
\notin Des(a,\prec)$ and $s \in Des(b,\prec)$. So consider
$a_{s-1}$. Our choice of $i$ and $s$ ensures that $|a_{s-1}|$ is
strictly greater than $i$. Thus by our definition of $\prec$, $i
\prec a_{s-1} \prec \overline{i}$. Hence $s-1 \in Des(a,\prec)$
and $s-1 \notin Des(b,\prec)$. Thus it again follows that
$\maj_{1,n-k+1}(b,\prec)-1= \maj_{1,n-k+1}(a,\prec)$
and hence
$$q^{pos(a)}q^{\maj_{1,n-k+1}(a,\prec)} = q^{pos(b)+1}q^{\maj_{1,n-k+1}(b,\prec)-1}
=q^{pos(b)}q^{\maj_{1,n-k+1}(b,\prec)}
$$
as desired.\\
\ \\
Thus our involution $I_k$ shows that
\begin{eqnarray}\label{AA21}
&&\sum_{a = a_1 \cdots a_n \in {\it A}_{\pm}^n}
q^{pos(a)}(-1)^{neg(a)} z_{|a_1|}\cdots z_{|a_n|}
q^{\maj_{1,n-k+1}(a,\prec)}t^{\rmaj_{n-k+1,n}(a,\prec)} = \\
&& \sum_{a = a_1 \cdots a_n \in {\it A}_{\pm}^n, I_k(a) =a}
q^{pos(a)}(-1)^{neg(a)} z_{|a_1|}\cdots z_{|a_n|}
q^{\maj_{1,n-k+1}(a,\prec)}t^{\rmaj_{n-k+1,n}(a,\prec)}.
\end{eqnarray}
But since the only words $a = a_1 \cdots a_n$ such that $I_k(a)
=a$ must have $|a_1|\cdots|a_{n-k+1}|$ be pairwise distinct, it
follows that the largest possible type of a monomial
$z_{|a_1|}\cdots z_{|a_n|}$ relative to the dominance order is
$(1^{n-k},k)$ since this is largest type of words with at least
$k+1$ distinct letters. Thus (\ref{AA1}) holds.

The proof of  (\ref{AA2}) is similar to the proof of (\ref{AA1}).
First we define  $\prec^*$ so that
$$
\overline{1} \prec^* \overline{2} \prec^* \cdots \prec^*
\overline{n} \prec^* n \prec^* n-1 \prec^* \cdots \prec^* 1.
$$
Then we define the weight $V(a)$ of a word $a = a_1 \cdots a_n \in
{\it A}_{\pm}^n$ by setting
$$V(a) = t^{pos(a)}(-1)^{neg(a)} z_{|a_1|}\cdots z_{|a_n|}
q^{\maj_{1,n-k+1}(a,\prec^*)}t^{\rmaj_{n-k+1,n}(a,\prec^*)}.$$
We define an involution $J_k$ on ${\it A}_{\pm}^n$ as follows.
Given a word $a = a_1 \cdots a_n \in {\it A}_{\pm}^n$, we look for
the smallest letter $j$ which is repeated in $|a_{n-k+1}| \cdots
|a_{n}|$. If there is no such letter, we let $J_k(a) =a$.
Otherwise, let $a_t$ be the left most occurrence of either $j$ or
$\overline{j}$ in $a_{n-k+1} \cdots a_{n}$. Now let $i$ be the
smallest letter in $|a_{t+1}| \cdots |a_{n}|$ and let $s$ be the
right most occurrence of either $i$ or $\overline{i}$ in $a_{t+1}
\ldots a_n$.

Thus $n-k+1 <t+1 \leq  s$ since $n-k+1 \leq t$. Then $J_k(a) = b= b_1
\cdots b_n$ where (i) $b_r = a_r$ if $r \neq s$, (ii) $b_s = i$ if
$a_s = \overline{i}$,and (iii) $b_s = \overline{i}$ if $a_s = i$.
In other words, $J_k(a)$ is the  result of changing $a_s$ to
$\overline{i}$ if $a_s =i$ or changing $a_s$ to $i$ if $a_s =
\overline{i}$. Clearly $J_k^2$ is the identity. We claim that
$V(a) = - V(b)$. Clearly $(-1)^{neg(a)} = -(-1)^{neg(b)}$ since we
changed the sign of one letter. Again we did not change the
absolute value of any letter so that $z_{|a_1|}\cdots z_{|a_n|} =
z_{|b_1|}\cdots z_{|b_n|}$. We did not change any of the letters
among $a_1 \cdots a_{n-k+1}$ so that $ q^{\maj_{1,n-k+1}(a,\prec^*)}
= q^{\maj_{1,n-k+1}(b,\prec^*)}$. Thus we need only show
that
\begin{equation}
t^{pos(a)}t^{\rmaj_{n-k+1,n}(a,\prec^*)} =
t^{pos(b)}t^{\rmaj_{n-k+1,n}(b,\prec^*)}
\end{equation}
There is no loss of generality in assuming that $a_s=i$.
Then there are two cases.\\
{\bf Case 1.} $s =n$. Since $i$ is the smallest element in
$|a_{t}| \cdots |a_{n}|$, it follows $a_{n-1} \preceq^* i$. We
claim that $n-1 \notin Des(a,\prec^*)$ and that $n-1 \in
Des(b,\prec^*)$.  Now if $|a_{n-1}| >i$, then clearly
$\overline{i} \prec^* a_{n-1} \prec^* i$ so that our claim holds.
If $a_{n-1} = i$, then $n-1 \notin Des(a,\prec^*)$ because two
equal positive letters do not cause a descent by our definitions.
However, $\overline{i} \prec^* i$ so that $n-1 \in
Des(b,\prec^*)$. Finally if $a_{n-1} = \overline{i}$, then $n-1
\in Des(b,\prec^*)$ since two equal negative letters cause a
descent while $ n-1 \notin Des(a,\prec^*)$ because  $\overline{i}
\prec^* i$. Thus it follows that
$\rmaj_{n-k+1,n}(b,\prec^*)-1=
\rmaj_{n-k+1,n}(a,\prec^*)$ and hence
$$t^{pos(a)}c^{\rmaj_{n-k+1,n}(a,\prec^*)} = t^{pos(b)+1}t^{\rmaj_{n-k+1,n}(b,\prec^*)-1}
=t^{pos(b)}t^{\rmaj_{n-k+1,n}(b,\prec^*)}
$$
as desired.\\
\ \\
{\bf Case 2.} $s <n$.  We can argue exactly as in Case 1 that $s-1
\notin Des(a,\prec^*)$ and $s-1 \in Des(b,\prec^*)$. So consider
$a_{s+1}$. Our choice of $i$ and $s$ ensures that $|a_{s+1}|$ is
strictly greater than $i$. Thus by our definition of $\prec^*$,
$\overline{i} \prec^* a_{s+1} \prec^* i$. Hence $s \in
Des(a,\prec^*)$ and $s \notin Des(b,\prec^*)$. Thus the only
places were $Des(a,\prec^*)$ and $Des(a,\prec^*)$ differ is on the
set $\{s-1,s\}$. We have $Des(a,\prec^*) \cap \{s-1,s\} = \{s\}$
and $Des(b,\prec^*) \cap \{s-1,s\} = \{s-1\}$. Since $n-(s-1) =
(n-s) +1$, it follows that $\rmaj_{n-k+1,n}(b,\prec^*)-1=
\rmaj_{n-k+1,n}(a,\prec^*)$ and hence
$$t^{pos(a)}t^{\rmaj_{n-k+1,n}(a,\prec^*)} = t^{pos(b)+1}t^{\rmaj_{n-k+1,n}(b,\prec^*)-1}
=t^{pos(b)}t^{\rmaj_{n-k+1,n}(b,\prec^*)}
$$
as desired.\\
\ \\

Thus our involution $J_k$ shows that
\begin{eqnarray}\label{AA11}
&&\sum_{a = a_1 \cdots a_n \in {\it A}_{\pm}^n}
t^{pos(a)}(-1)^{neg(a)} z_{|a_1|}\cdots z_{|a_n|}
q^{\maj_{1,n-k+1}(a,\prec^*)}t^{\rmaj_{n-k+1,n}(a,\prec^*)} = \\
&& \sum_{a = a_1 \cdots a_n \in {\it A}_{\pm}^n, J_k(a) =a}
t^{pos(a)}(-1)^{neg(a)} z_{|a_1|}\cdots z_{|a_n|}
q^{\maj_{1,n-k+1}(a,\prec^*)}t^{\rmaj_{n-k+1,n}(a,\prec^*)}.
\end{eqnarray}
But since the only words $a = a_1 \cdots a_n$ such that $J_k(a)
=a$ must have $|a_{n-k+1}|\cdots|a_{n}|$ be pairwise distinct, it
follows that the largest possible type of a monomial
$z_{|a_1|}\cdots z_{|a_n|}$ relative to the dominance order is
$(1^{k-1},n-k+1)$ since this is largest type of words with at least
$k$ distinct letters. Thus (\ref{AA2}) holds.

\qed

\section{Final Remarks}\label{QS}

\subsection{Haglund Statistics}\label{QS1}

\begin{figure}[tbp]
\centering
$$
\epsfysize=.7in \epsffile{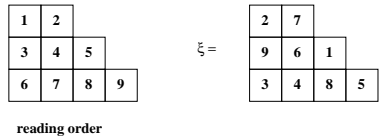}
$$
\caption{The reading order and a filling of
$(4,3,2)$.}\label{fig:mac3}\end{figure}

Let $\xi$ be a filling of the Ferrers diagram of a partition $\mu$
with the numbers $1, \ldots, n$. For any cell $u = (i,j)\in
F_\mu$, let $\xi(u)$ be the entry in cell $u$. We say that $u =
(i,j)\in F_\mu$ is a {\em descent} of $\xi$, written $u\in
Des(\xi)$, if $i > 1$ and $\xi((i,j)) \geq \xi((i-1,j))$. Then
$maj(\xi) = \sum_{u \in Des(\xi)} (leg(u) +1)$. Two cells $u,v \in
F_\mu$ {\em attack} each other if either
\begin{itemize}
\item[(a)] they are in the same row, i.e. $u=(i,j)$ and $v=(i,k)$, or
\item[(b)] they are in consecutive rows, with the cell in the
upper row strictly to the right of the one in the lower row, i.e.
$u=(i+1,k)$ and $v=(i,j)$ where  $j < k$.
\end{itemize}
The {\em reading order} is the total ordering on the cells of
$F_\mu$ given by reading the cells row by row from top to bottom,
and left to right within each row. For example, the reading order
of $(2,3,4)$ is depicted on the left in Figure~\ref{fig:mac3}. An
{\em inversion} of $\xi$ is a pair of entries $\xi(u) > \xi(v)$
where $u$ and $v$ attack each other and $u$ precedes $v$ in the
reading order. We then define $Inv(\xi)= \{ \{u,v\}: \xi(u) >
\xi(v) \ \mbox{is an inversion}\}$ and $inv(\xi) = |Inv(\xi)| -
\sum_{u \in Des(\xi)} arm(u)$.

For example, if $\xi$ is the filling of shape $(2,3,4)$ depicted
in Figure~\ref{fig:mac3}, then $Des(\xi) = \{(2,1),(2,2),(3,2)\}$.
There are four inversion pairs of type (a), namely
$\{(2,1),(2,2)\}$, $\{(2,1),(2,3)\}$, $\{(2,2),(2,3)\}$, and
$\{(1,3),(1,4)\}$, and one inversion pair of type (b), namely
$\{(2,2),(1,1)\}$. Then one can check that $|Inv(\xi)| = 5$,
$maj(\xi) = 5$ and $inv(\xi) = 2$. Finally, we can identify $\xi$
with a permutation by reading the entries in the reading order. In
the example of Figure~\ref{fig:mac3}, $\xi = 2~7~9~6~1~3~4~8~5$.
Then we let $D(\xi) = Des(\xi^{-1})$. In our example, $\xi^{-1} =
5~1~6~7~9~4~2~8~3$ so that $D(\xi) = \{1,5,6,8\}$.

Recently, Haglund, Haiman and Loehr~\cite{HHL05a,HHL05b} proved
Haglund's conjectured combinatorial interpretation \cite{Hag04} of
${\tilde H}_{\mu}(\ox;q,t)$ in terms of quasi-symmetric functions.
That is, given a non-negative integer $n$ and a subset $D
\subseteq \{1, \ldots, n-1\}$, Gessel's quasi-symmetric function
of degree $n$ in variables $x_1, x_2, \ldots $ is defined by the
formula
\begin{equation}\label{quasi}
Q_{n,D}(\ox): = \sum_{\overset{a_1\leq a_2 \leq \cdots \leq
a_n}{a_i=a_{i+1} \Rightarrow i \notin D}} x_{a_1} x_{a_2} \cdots
x_{a_n}.
\end{equation}
Then Haglund, Haiman and Loehr \cite{HHL05b} proved
\begin{equation}\label{quasiHmu}
{\tilde H}_{\mu}(\ox;q,t) = \sum_{\xi:\mu \simeq \{1,\ldots, n\}}
q^{inv(\xi)}t^{maj(\xi)} Q_{n,D(\xi)}(\ox).
\end{equation}
Here the sum runs over all fillings $\xi$ of the Ferrers diagram
of $\mu$ with the numbers $1, \ldots, n$.

\subsection{Relations with the Combinatorial Interpretation of Macdonald Polynomials
}\label{QS2}

The Hilbert series of ${\bf H}_\mu$ is equal to the coefficient of
$x_1x_2 \cdots x_n$ in ${\tilde H}_{\mu}(\ox;q,t)$. Since the
coefficient of $x_1 x_2 \cdots x_n$ in any quasi-symmetric
function $Q_{n,D}(\ox)$ is $1$, it follows that the Hilbert series
of ${\bf H}_\mu$ is given by
$$
\sum_{k,r} \dim {\bf H}_\mu^{(h,k)} q^h t^k
= {\tilde H}_{\mu}(\ox;q,t)|_{x_1x_2 \cdots x_n} \\
= \sum_{\xi:\mu \simeq \{1,\ldots, n\}} q^{inv(\xi)}t^{maj(\xi)} ,
$$
where the sum runs over all fillings $\xi$ of the Ferrers diagram
of $\mu$ with the numbers $1, \ldots, n$. No known basis realizes
this remarkable identity for general ${\bf H}_\mu$. The $k$-th
Haglund basis described in Subsection~\ref{s.k-AH} above provides
such a basis when $\mu$ is of hook shape.


\begin{thebibliography}{99}

\bibitem{TAMS}
R.\ M.\ Adin, F.\ Brenti and Y.\ Roichman, {\rm Descent
representations and multivariate statistics}, {\em Trans.\ Amer.\
Math.\ Soc.}~{\bf 357} (2005), 3051--3082.

\bibitem{APR1}
R.\ M.\ Adin, A.\ Postnikov and Y.\ Roichman, {\rm Hecke algebra
actions on polynomial rings}, {\em J.\ Algebra}~{\bf 233} (2000),
594--613.

\bibitem{APR2}
R.\ M.\ Adin, A.\ Postnikov and Y.\ Roichman, {\rm On characters
of Weyl groups}. {\em Discrete Math.}~{\bf 226} (2001), 355-358.

\bibitem{Allen1}
E.\ E.\ Allen, {\rm The descent monomials and a basis for the
diagonally symmetric polynomials}, {\em J.\ Alg.\ Combin.}~{\bf 3}
(1994), 5--16.

\bibitem{Allen2}
E.\ E.\ Allen, {\rm Bitableaux bases for the diagonally invariant
polynomial quotient rings}, {\em Adv. Math.}~{\bf 130} (1997),
242--260.

\bibitem{Aval}
J.-C. Aval, {\rm Monomial bases related to the $n!$ conjecture},
Discrete Mathematics~{\bf 224}, (2000), 15-35.

\bibitem{Ba}
H.\ Barcelo, {\rm Young straightening in a quotient $S_n$-module},
{\em J.\ Alg.\ Combin.}~{\bf 2} (1993), 5--23.

\bibitem{BReg87}
A.\ Berele and A.\ Regev, {\rm Hook Young diagrams and their
applications to combinatorics and to representations of Lie
superalgebras}, {\em Adv. in Math.}~{\bf 64} (1987), 118--175.

\bibitem{BRem85}
A.\ Berele and J.\ B.\ Remmel, {\rm Hook flag characters and their
combinatorics}, {\em J.\ Pure Appl.\ Alg.}~{\bf 35} (1985),
225--245.

\bibitem{BBGHT}
N.\ Bergeron, F.\ Bergeron, A.\ M.\ Garsia, M.\ Haiman, and G.\
Tesler, {\rm Lattice diagram polynomials and extended Pieri
rules}, {\em Adv. Math.}~{\bf 142} (1999), 244--334.

\bibitem{GH93}
A.\ M.\ Garsia and M.\ Haiman, {\rm A graded repesentation model
for the Macdonald polynomials}, {\em Proc.\ Nat.\ Acad.\
Sci.}~{\bf 90} (1993), 3607--3610.

\bibitem{GH95}
A.\ M.\ Garsia and M.\ Haiman, {\rm Factorizations of Pieri rules
for Macdonald polynomials}, {\em Discrete Math.}~{\bf 139} (1995),
219--256.

\bibitem{GH1}  A.\ M.\ Garsia and M.\ Haiman, {\rm Some
natural bigraded $S_n$-modules}, {\em Electronic Journal of
Combinatorics}~{\bf 3} (1996), R24.

\bibitem{GHOH} A.\ M.\ Garsia and M.\ Haiman,
{\rm Orbit Harmonics and Graded Representations}, Research
Monograph to appear as part of the collection published by the
Lab. de. Comb. et Informatique Math\'ematique, edited by S. Brlek,
U. du Qu\'ebec \'a Montr\'eal.

\bibitem{GP}
A.\ M.\ Garsia and C.\ Procesi, {\rm On certain graded
$S_n$-modules and the $q$-Kostka polynomials}, {\em Adv.\
Math.}~{\bf94} (1992), 82--138.

\bibitem{GR}
A.\ M.\ Garsia and J.\ Remmel, {\rm Shuffles of permutations and
the Kronecker product}, {\em Graphs and Combinatorics}~{\bf 1}
(1985), 217--263.

\bibitem{GS}
A.\ M.\ Garsia and D.\ Stanton, {\rm Group actions of
Stanley-Reisner rings and invariants of permutation groups}, {\em
Adv.\ Math.}~{\bf 51} (1984), 107--201.

\bibitem{Ge2}
I.\ M.\ Gessel, {\rm Multipartite $P$-partitions and inner
products of Schur functions}, {\em Contemp. Math.}~{\bf 34}
(1984), 289--302.

\bibitem{Hag04}
J.\ Haglund, {\rm A Combinatorial Model for the Macdonald
polynomials}, {\em Proc.\ Nat.\ Acad.\ Sci.}~{\bf 101} (2004),
16127--16131.

\bibitem{HHL05a}
J.\ Haglund, M.\ Haiman, and N.\ Loehr, {\rm Combinatorial theory
of Macdonald polynomials I: Proof of Haglund's formula}, {\em
Proc.\ Nat.\ Acad.\ Sci.}~{\bf 102} (2005), 2690--2696.

\bibitem{HHL05b}
J.\ Haglund, M.\ Haiman, and N.\ Loehr, {\rm A Combinatorial
formula for the Macdonald  polynomials}, {\em J.\ Amer.\ Math.\
Soc.}~{\bf 18} (2005), 735--761.

\bibitem{Hai99}
M.\ Haiman, {\rm Macdonald polynomials and geometry}, in: {\em New
Perspectives in algebraic geometry (Berkeley, CA, 1996-97)}, {\em
Math.\ Sci.\ Res.\ Inst.\ Publ.}, vol.~38, pp.\ 207--254,
Cambridge Univ. Press, Cambridge (1999).

\bibitem{Hai01}
M.\ Haiman, {\rm Hilbert schemes, polygraphs, and the Macdonald
positivity conjecture}, {\em J.\ Amer.\ Math.\ Soc.}~{\bf 14}
(2001), 941--1006.

\bibitem{Hum}
J.\ E.\ Humphreys, {\rm Reflection Groups and Coxeter Groups},
{\em Cambridge Studies in Advanced Mathematics}~{\bf 29},
Cambridge Univ.\ Press, Cambridge, 1990.

\bibitem{KL}
D.\ Kazhdan and G.\ Lusztig, {\rm Representations of Coxeter
groups and Hecke algebras}, {\em Invent.\ Math.}~{\bf 53} (1979),
165--184.

\bibitem{Knuth}
D.\ E.\ Knuth, {\rm Permutations, matrices and generalized Young
tableaux}, {\em Pacific J.\ Math.}~{\bf 34} (1971), 709--727.

\bibitem{Mac88}
I.\ G.\ Macdonald, {\rm A new class of symmetric functions}, {\em
Publ.\ I.R.M.A.\ Strasbourg, Actes $20^e$ S\'{e}minaire
Lotharingien}~(1988), 131--171.

\bibitem{Mac95}
I.\ G.\ Macdonald, {\rm Symmetric Functions and Hall Polynomials},
$2^{nd}$ ed., Oxford Univ.\ Press (1995).

\bibitem{RS02}
A.\ Regev and T.\ Seeman, {\rm Shuffle invariance of the super-RSK
algorithm}, {\em Adv.\ Appl.\ Math.}~{\bf 28} (2002), 59--81.

\bibitem{Rem84}
J.\ B.\ Remmel, {\rm The combinatorics of $(k,l)$-hook Schur
functions}, {\em Contemp.\ Math.\ AMS}, vol.~34 (Comb.\ and Alg.)
(1984), pp.\ 253--287.

\bibitem{Rem90}
J.\ B.\ Remmel, {\rm A bijective proof of a factorization theorem
for $(k,l)$-hook Schur functions}, {\em Lin.\ and Multilin.\
Alg.}~{\bf 28} (1990), 253--287.

\bibitem{Rob38}
G.\ de B.\ Robinson, {\rm On the representations of $S_n$}, {\em
Amer.\ J.\ Math.}~{\bf 60} (1938), 745--760.

\bibitem{Sch61}
C.\ E.\ Schensted, {\rm Longest increasing and decreasing
subsequences}, {\em Canad.\ J.\ Math.}~{\bf 13} (1961), 179--191.

\bibitem{Sag91}
B.\ Sagan, {\em The Symmetric Group: Representations,
Combinatorial Algorithms, and Symmetric Functions}, Wadsworth and
Brooks/Cole, 1991.

\bibitem{So}
L. Solomon, {\rm The orders of the finite Chevalley groups}, {\em
J. Algebra}~{\bf 3} (1966), 376--393.

\bibitem{St79}
R.\ P.\ Stanley, {\rm Invariants of finite groups and their
applications to combinatorics}, {\em Bull.\ Amer.\ Math.\ Soc.\
(new series)}~{\bf 1} (1979), 475--511.

\bibitem{St82}
R.P. Stanley, {\rm Some aspects of group acting on finite posets},
{\em J. Combin. Theory Ser. A}~{\bf 32} (1982), 132--161.

\bibitem{Ste94}
J.\ R.\ Stembridge, {\rm Some particular entries of the
two-pararemeter Kostka matrix}, {\em Proc.\ Amer.\ Math.\
Soc.}~{\bf 121} (1994), 469--490.

\bibitem{Stg}
R.\ Steinberg, {\rm On a theorem of Pittie}, {\em Topology}~{\bf
14} (1975), 173--177.

\end{thebibliography}
\end{document}